\documentclass[12pt]{amsart}
\usepackage{graphicx} % Required for inserting images
\usepackage{amsmath,latexsym}
\usepackage{todonotes}
\usepackage{mathabx}
\usepackage{mathtools}
\usepackage{bm}
\usepackage{enumitem}

\DeclarePairedDelimiter{\ceil}{\lceil}{\rceil}

\newtheorem{proposition}{Proposition}
\newtheorem{theorem}{Theorem}
\newtheorem{lemma}{Lemma}
\theoremstyle{definition}
\newtheorem{assumption}{Assumption}

\title[Enhanced dissipation]{Enhanced dissipation by advection and applications to PDEs}
\author[]{Anna L. Mazzucato}
\address{Department of Mathematics, Penn State University, University Park, PA
16802, USA}
\email{alm24@psu.edu}
\author[Y. Feng {\em et Al.}]{Yuanyuan Feng}
\address{School of Mathematical Sciences,  Key Laboratory of MEA (Ministry of Education) \& Shanghai Key Laboratory of PMMP,  East China Normal University, Shanghai 200241, China }
\email{yyfeng@math.ecnu.edu.cn}
\author[]{Camilla Nobili}
\address{ School of Mathematics and Physics, University of Surrey, Guildford, GU2 7XH, UK }
\email{c.nobili@surrey.ac.uk}

\usepackage[backend=bibtex,sortcites=true,style=numeric-comp]{biblatex}
\usepackage[left=3cm,right=3cm,top=2cm,bottom=2cm]{geometry}

\date{\today}

\begin{document}
%\tableofcontent

\begin{abstract}
This survey provides a concise yet comprehensive overview on enhanced dissipation phenomena, transitioning seamlessly from the physical principles underlying the interplay between advection and diffusion to their rigorous mathematical formulation and analysis. The discussion begins with the standard theory of enhanced dissipation, highlighting key mechanisms and results, and progresses to its applications in notable nonlinear PDEs such as the Cahn-Hilliard and Kuramoto-Sivashinsky equations.
\end{abstract}

\maketitle

\section{Introduction} \label{s:intro}

This survey concerns the phenomenon of {\em enhanced dissipation} by the addition of linear advection to dissipative systems, some of its consequences, and applications to non-linear partial differential equations arising in fluid mechanics and related areas. 
It showcases results obtained by the authors and their co-authors, along with a review of key results in the literature.

Enhanced dissipation means that dissipation acts more efficiently in the presence of transport by a given background flow, provided the flow satisfies some conditions. Dissipation can be modeled by rather general operators, but for simplicity in this work we will refer only to the Laplace or bi-harmonic operators. Enhanced dissipation can be measured in terms of the so-called {\em dissipation time}, informally the smallest time on which the advection-(hyper)diffusion dynamical system reduces the size of the solution in half, or in terms of decay rates in time for the solution operator of the advection-(hyper)diffusion equation. Then we speak of enhanced dissipation when the dissipation time of the advection-(hyper)diffusion system is strictly shorter than that of (hyper)diffusion alone.
The concept of enhance dissipation is related to other physically observed and relevant phenomena as {\em Taylor dispersion} and {\em mixing} by the background flow.

While mixing and enhanced dissipation have been extensively studied, for instance in the context of ergodic theory and probability (see e.g. \cite{DaPratoBook}), there has been recently a renewed interest in these important phenomena in the context of partial differential equations, transport theory, and fluid mechanics, in connection with scalar turbulence and also stability near steady flows for instance (we refer to the survey \cite{CZCIMSurvey} and the lecture notes \cite{LevicoNotes} for more details on these aspects), starting from the seminal contribution of Constantin {\em et Al.} \cite{CKRZ08}, and of Doering {\em et Al.} \cite{DoeringThiffeault06}. In particular, it was recognized that enhanced dissipation plays a fundamental role in regulating the behavior of solutions to certain dissipative systems, modeled by non-linear partial differential equations. For instance, in aggregation models, such as the Keller-Segel-Patlak chemotaxis model, addition of a convective term by a mixing flow can suppress finite-time blow-up of solutions with large total mass \cites{KiselevXu16,IyerXuEA21,BH18,HK22,HKY23ArXiV}. It can suppress phase separation in phase-field models such as those described by the Allen-Cahn or Cahn-Hilliard equation \cite{FengFengEA20}, and prolong the life of solutions in dissipation systems with long-wave instabilities, such as in the Kuramoto-Sivashinsky equation in two and higher space dimensions \cites{FM22,CZDFM21,FengShiEA22}, or lead to symmetrization in the two-dimensional Navier-Stokes equations at long times \cite{CZEW20} (see e.g. \cites{DL22,BMV19,GNRS20,WZZ20,MZ20,J23,CL24,DZ23,W25,Z21} also for connections with stability near steady flows and the phenomenon of inviscid damping, and \cite{JW22} for dissipation enhancement by vortex strechting).

There is also an important link with the phenomenon of {\em anomalous dissipation}, that is, the observed phenomenon that the energy dissipation rate of the passive scalar seems to reach a plateau at a non-zero value, as diffusivity tends to zero, the analog in scalar turbulence of the so-called ``zeroth-law of turbulence" in  Kolmogorov's theory. We refer the reader to recent rigorous advances in these areas \cites{EL24,BCCDLS24,HL23ArXiv,CCS23,R24,BDL23,K24ArXiV,DEIJ22}, which primarily use ideas from mixing and irregular transport, \cite{AV25,BSLW23ArXiV}, which use ideas from homogenization, and references therein (see also \cite{DG12} from a turbulence-inspired approach). The limit of zero viscosity or zero diffusivity has been investigated  in this context also as a selection principle for weak solutions (among the extensive literature, we mention  \cites{BTW12,BLNNT13,HL23ArXiv,BCC22,CCS23} and references therein).
Since we work with positive diffusivity, to keep this introduction brief we do not discuss the fundamental related problem of energy dissipation at zero diffusivity or viscosity, that is, whether weak solutions of the linear transport or Euler equations dissipate energy (in the literature, this problem is often called anomalous dissipation, while the zeroth law of turbulence is referred to as {\em dissipation anomaly}) .

The survey is organized as follows. Section \ref{s:enhanced} deals with enhanced dissipation for advection-(hyper)diffusion equations. In particular, Section \ref{s:dispersion} introduces the concept of enhanced dissipation and Taylor's dispersion, especially for shear flows,  and compares them, while Section \ref{s:resolvent} discusses  methods, more specifically {\em hypocoercivity} and {\em resolvent estimates}, to quantify enhanced dissipation. Applications to non-linear dissipative equations are given in Section \ref{s:applications}, with Section \ref{s:CHE} dealing with enhanced dissipation by mixing flows, Cahn-Hilliard and aggregation models, while Section \ref{s:KSE} presents results concerning global existence for the Kuramoto-Sivashinksy equation with advection in two space dimensions.

\subsection*{Acknowledgments} The authors thank Michele Coti Zelati and Gautam Iyer for fruitful discussions on the results surveyed in this work. Y. Feng was partially supported by National Key Research and Development Program of China (2022YFA1004401), NSFC 12301283, Shanghai Sailing program 23YF1410300, Science and Technology Commission of Shanghai Municipality (22DZ2229014). A. M. was partially supported by the US National Science Foundation Grants DMS-1909103, DMS-2206453, and Simons Foundation Grant 1036502.

\section{Enhanced dissipation for advection-diffusion equations} \label{s:enhanced}

\subsection{Enhanced dissipation and Taylor dispersion} \label{s:dispersion}

In this section, we consider the advection-diffusion equation for a scalar tracer $\theta = \theta(\bm{x}, t)$:
\begin{equation}\label{ad-eq}
\partial_t \theta + \bm{u} \cdot \nabla \theta = \kappa \Delta \theta,
\end{equation}
in simple geometries such as the torus, the periodic channel, and the infinite channel. Here, $\kappa$ is the molecular diffusivity, and the spatial variable $\bm{x}$ is $\bm{x} = (x_1, x_2, x_3)$ in three dimensions or $\bm{x} = (x_1, x_2)$ in two dimensions.

We begin by clarifying the distinction between diffusion and dispersion, both of which are crucial in understanding mixing and enhanced dissipation. \textit{Diffusion} (or molecular diffusion) is a microscopic process resulting from the random motion of individual molecules moving from regions of higher concentration to regions of lower concentration. In contrast, \textit{dispersion} is a macroscopic phenomenon characterized by the spreading of an initially localized concentration profile due to a combination of advective velocity gradients, molecular diffusion, and geometric boundary effects.

In what follows, we will provide a mathematical explanation of \emph{shear dispersion} (also known as Taylor’s dispersion), following the seminal work of Taylor~\cite{Taylor53}, and discuss its role in enhanced dissipation.

\subsubsection{Shear dispersion and diffusivity enhancement}\label{enhanced-dispersion}
In his seminal work in \cite{Taylor53} Taylor developed a quantitative theory for the spread of tracer along the axis of a a pipe with radius $a$, described, in Cartesian coordinates, by \eqref{ad-eq} with  axial velocity along the pipe $$\bm{u}=(u(x_2,x_3), 0, 0) \quad \mbox{and} \quad  u(x_2, x_3)= 2\frac{\bar u}{a^2}(a^2-x_2^2-x_3^2),$$
where $\bar{u}$ is the mean velocity of the flow.
\begin{figure}
\centering
\includegraphics[scale=0.5]{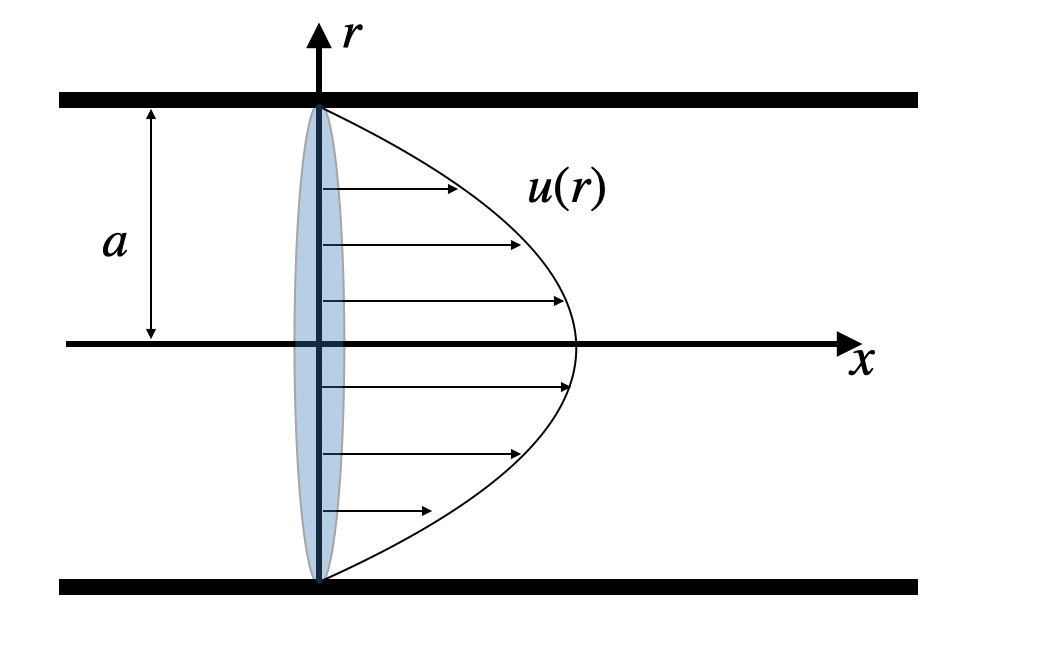}
\caption{Illustration of the channel setting and the velocity field considered by Taylor in \cite{Taylor53}}
\end{figure}
In the more convenient cylindrical coordinates $(r,\phi, x_3)$, the velocity reads  $(u(r), 0, 0)$,
where 
\begin{equation}\label{velocity}
u(r)=2 \bar{u} (1-\frac{r^2}{a^2})
\end{equation}
and $\bar{u}=\frac{2}{ a^2}\int_0^a  u \,r\, dr $. 
The evolution for $\theta(x_1, r, t)$ is described by 
\begin{equation}\label{ad-cc}
\partial_t \theta+u(r)\partial_{x_1} \theta=\kappa (\partial_{x_1}^2\theta+\frac{1}{r}\partial_r(r\partial_r\theta))\,,
\end{equation}
where the term depending on the angular variable, $\frac{1}{r^2}\partial_{\phi}^2\theta$, was neglected under the assumption of an axisymmetric flow.
We decompose $\theta$ in a cross-sectional average and $r$-dependent part, i.e.,
$$\theta(x_1 ,r ,t)=\bar{\theta}(x_1, t)+\theta'(x_1, r, t) \quad \mbox{ where }\quad \bar{\theta}=\frac{2}{a^2}\int_0^a \theta\, r \, dr \mbox{ and } \bar{\theta'}=0\,,$$
and, in analogy, the velocity field is decomposed as
$$u(x_1, r)=\bar{u}(x_1)+u'(x_1, r). $$% \quad \mbox{ where }\quad \bar{u}=\frac{2}{a^2}\int_0^a u\, r \, dr \mbox{ and } \bar{u'}=0\,.$$
Moreover we assume the boundary condition 
\begin{equation}\label{bctheta'}
    \partial_{r}\theta=0 \mbox{ at } r=a\,.
\end{equation}
    
%By imposing the decompositions into \eqref{ad-eq} it is easy to see that 
%and $\theta'$ solves 
%\begin{equation}\label{eqtheta'}
%\partial_t \theta'+u'\partial_{x_3}\bar{\theta}+u\partial_{x_3}\theta'- \widebar{u'\partial_{x_3}\theta'}=\kappa \Delta \theta'=\kappa (\frac{1}{r}\partial_r(r\partial_r \theta')+\partial_{x_3}^2 \theta')\,.
%\end{equation}
Imposing this decomposition in \eqref{ad-cc} gives
\begin{equation}
\partial_t(\bar{\theta}+\theta')+(\bar{u}+u')\partial_{x_1}(\bar{\theta}+\theta')=\kappa (\partial_{x_1}^2(\bar{\theta}+\theta')+\frac{1}{r}\partial_r(r\partial_r\theta')),
\end{equation}
where we used that $\partial_r \bar{\theta}=0$. Now, taking the average in $r$ of the above equation yield
\begin{equation}\label{eqthetabar}
\partial_t \bar{\theta}+\bar{u}\partial_{x_1}\bar{\theta}+\widebar{ u' \partial_{x_1}\theta'}=\kappa \partial_{x_1}^2 \bar{\theta}.
\end{equation}
Following Taylor's arguments we will show how the mixing term $\widebar{ u' \partial_{x_1}\theta'}$ enhances dispersion. Specifically, by using some physical assumptions we will show that this equation reduces to
\begin{equation*}
\partial_t \bar{\theta} + \bar{u}\partial_{x_1} \bar{\theta} =\kappa_{\rm{eff}} \partial_{x_1}^2\bar{\theta}
\end{equation*}
which is a one-dimensional advection-diffusion equation with a new \textit{effective diffusivity} which reflects the enhanced dispersion phenomena due to shear.

We view the flow as an observer that is moving at the mean velocity and thus make the following change of variable 
$$\xi=x_1-\bar{u}t, \qquad \tau=t,$$
in which we notice that the origin moves with the mean flow velocity.
Our equation then transforms to 
\begin{equation*}
\partial_{\tau}(\bar{\theta}+\theta')+u'\partial_{\xi}(\bar{\theta}+\theta')=\kappa (\frac{1}{r}\partial_r(r\partial_r \theta'))\,.
\end{equation*}
We notice that in the new coordinates system $u'$ is the only observable velocity.
We can neglect longitudinal diffusion because the rate of spreading along the flow direction due to velocity differences greatly exceeds that due to molecular diffusion. So we impose 
\begin{equation*}
    u'\partial_{\xi}(\bar{\theta}+\theta')\gg \kappa \partial_{\xi}^2(\bar{\theta}+\theta')
\end{equation*}
and we are left to analyze
\begin{equation}\label{dec2}
\partial_{\tau}(\bar{\theta}+\theta')+u'\partial_{\xi}(\bar{\theta}+\theta')=\kappa (\frac{1}{r}\partial_r(r\partial_r \theta'))\,.
\end{equation}
Averaging \eqref{dec2} in the radial direction and using  $\overline{r^{-1}\partial_r(r\partial_r\theta')}=0$, $\overline{\theta'}=0$, 
 and $\overline{u'}=0$, we obtain thanks to \eqref{bctheta'} that
$$\partial_{\tau}\bar{\theta}+\widebar{u'\partial_{\xi}\theta'}=0\,.$$
Subtracting this equation from \eqref{dec2} we get 
\begin{equation*}
\partial_{\tau}\theta'+u'\partial_{\xi}\bar{\theta}-\overline{u'\partial_{\xi}\theta'}+u'\partial_{\xi}\theta'=\kappa (\frac{1}{r}\partial_r(r\partial_r \theta'))\,.
\end{equation*}
We now assume that $\bar{\theta}, \theta'$ are well behaved, slowly varying and 
\begin{equation}\label{ass2}
    \bar{\theta}\gg \theta'\,,
\end{equation}
then
\begin{eqnarray}\label{approx}
u'\partial_{\xi}\bar{\theta}&\gg& u' \partial_{\xi}\theta'\\
u'\partial_{\xi}\bar{\theta}&\gg& \overline{u' \partial_{\xi}\theta'}\notag
\end{eqnarray}
and the equation further reduces to
\begin{equation}\label{crucial-reduction}
\partial_{\tau}\theta'+u'\partial_{\xi}\bar{\theta}=\kappa (\frac{1}{r}\partial_r(r\partial_r \theta'))\,.
\end{equation}
It is important to observe that for  \eqref{approx} to hold we need that the aspect ratio $\frac{a}{L}\ll 1$. 
%This can be verified by (a posteriori) scaling analysis.
The next crucial assumption made by Taylor is that, after sufficient time has elapsed for 
molecular diffusion in the transverse direction to fully develop, specifically for times satisfying
\[
t \gg \frac{a^2}{\kappa},
\]
the term $\partial_{\tau}\theta'$ can be neglected. Under this approximation, the governing equation reduces to
\[
u' \partial_{\xi}\bar{\theta} = \kappa \left( \frac{1}{r} \partial_r\left(r \partial_r \theta'\right) \right),
\]
representing a balance between longitudinal  transport and cross-sectional diffusive transport.

%In this regime we expect that the transverse differences in concentration represented by $\theta'$ are much less than the axial variations contained in $\bar{\theta}$ so that $\bar{\theta}\gg \theta'$. When $t\gg \frac{a^2}{\kappa}$ the gradients in the $z-$ direction become negligible with respect to the one in the $r-$direction, so equation \eqref{eqtheta'} essentially reduces to 
%$$(U(r)-\bar{U})\partial_{x_3}\bar{\theta}=\kappa (\frac{1}{r}\partial_r(r\partial_r \theta'))\,.$$
We can now finally solve this equation to find $\theta'$. Using \eqref{velocity} and integrating two times in $r$, we have 
\begin{equation}
    \theta'=\frac{\bar{u}\partial_{\xi}\bar{\theta}}{\kappa}\left(\frac{r^2}{4}-\frac{r^4}{8a^2}+A\log r+B\right)\,,
\end{equation}
where $A=0$ and $B=-\frac{1}{12}a^2$ are determined by the fact that $\theta'$ must be regular at $r=0$ and $\bar{\theta'}=0$.

We can now finally rewrite the term $\widebar{ u'\partial_{x_1}\theta'}$ in the original variables as 
$$\overline{ u' \partial_{x_1}\theta'}=-\frac{a^2 \bar{u}^2}{48 \kappa}\partial_{x_1}^2\bar{\theta}\,.$$
Inserting this expression back into \eqref{eqthetabar} we obtain 
$$\partial_t \bar{\theta}+\bar{u}\partial_{x_1}\bar{\theta}=\left(\kappa+\frac{a^2\bar{u}^2}{48 \kappa}\right)\partial_{x_1}^2 \bar{\theta}.$$
This computation show that  the large times behavior of the tracer is described by the one-dimensional advection-diffusion equation 
\begin{equation}\label{tay-eq}
\partial_t \bar{\theta} + \bar{u}\partial_{x_1} \bar{\theta} =\kappa_{\rm{eff}} \partial_{x_1}^2\bar{\theta}
\end{equation}
where $\kappa_{\rm{eff}}$ is the \textit{effective diffusivity} 
$$\kappa_{\rm{eff}}=\left(\kappa+\frac{a^2\bar{u}^2}{48 \kappa}\right)\,.$$

With this argument, Taylor obtained an enhancement of the diffusivity by a factor 
\[
\frac{a^2\bar{u}^2}{48\kappa},
\]
which is notably inversely proportional to the molecular diffusivity $\kappa$. Although the underlying idea is conceptually straightforward, this example provides crucial insights into the mechanisms responsible for diffusion enhancement, and it will be revisited later in this review article. We observe that when $\kappa$ is small, the effective diffusivity becomes extremely large. This scenario can be referred to as a strongly diffusive regime (see also the illustrative example in Section~2 of \cite{rhines1983rapidly}).  

\begin{figure}
\centering
\includegraphics[scale=0.5]{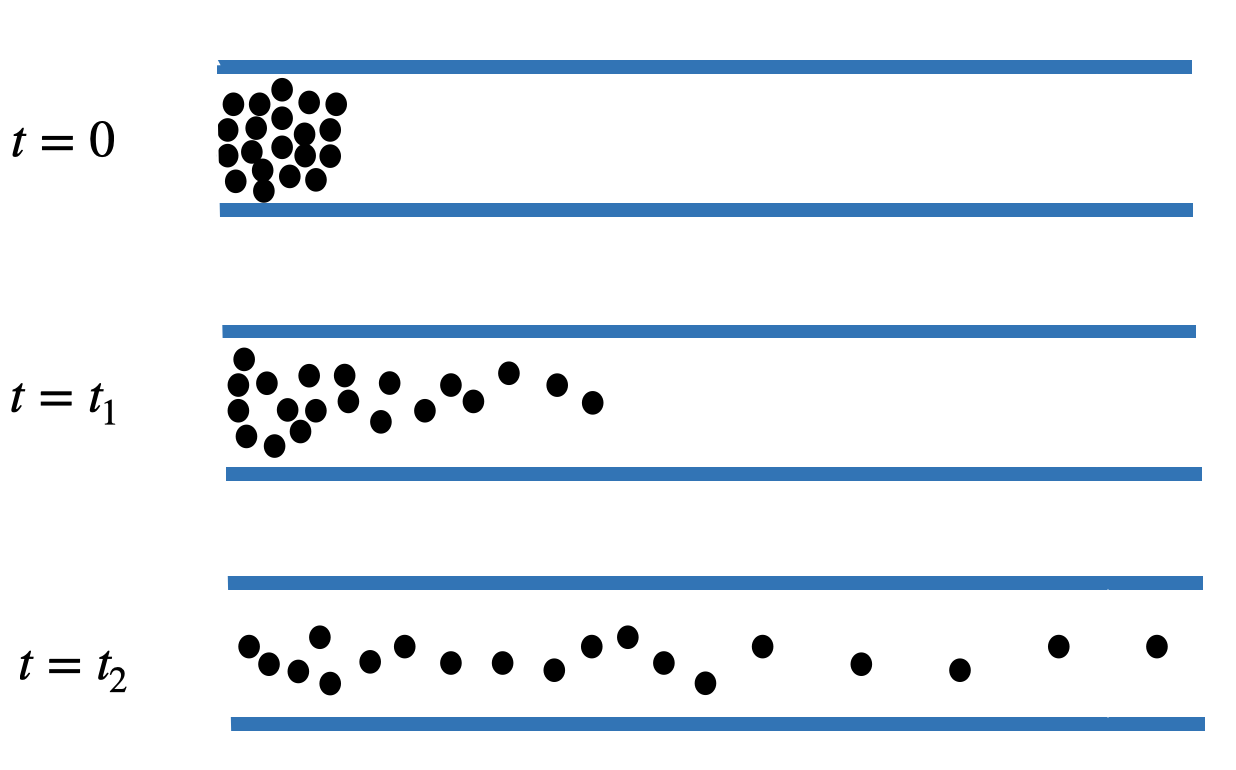}
\caption{Illustration of dispersion phenomena}
\end{figure}

We now turn to an example from \cite{rhines1983rapidly}, where the interplay of velocity gradients and weak diffusivity leads to an accelerated spatial flux of the tracer $\theta$. This example elucidates the behavior of the system when $\kappa$ is small. We consider equation \eqref{ad-eq} on $\mathbb{T}^2$ with $\bm{u} = (\alpha x_2, 0)$ and initial condition  
\begin{equation}\label{id}
    \theta(x_1,x_2,0) = \cos(k x_1). 
\end{equation}
If $\kappa = 0$, the solution of \eqref{ad-eq} is 
\[
\theta(x_1, x_2, t) = \cos(k(x_1 - \alpha x_2 t)).
\]
By using the ansatz $\theta(x_1, x_2, t) = A(t)\cos(k(x_1 - \alpha x_2 t))$, it is straightforward to find the solution for $\kappa \neq 0$ (see \cites{rhines1983rapidly,doering2017lectures} for details):
\[
\theta(x_1,x_2,t) = \exp\left(-\kappa k^2 t - \tfrac{1}{3}\kappa k^2 \alpha^2 t^3\right)\cos(k(x_1 - \alpha x_2 t)).
\]
From this explicit solution, we immediately infer that the level sets of the initial data become tilted. Setting $\alpha=0$ recovers the  solution  of the standard heat equation with the given initial condition \eqref{id}.

Next, we estimate the $L^2$-norms of the solution and its gradient:
\[
\iint_{\mathbb{T}^2} |\theta|^2\, d\bm{x} \leq (2\pi)^2\exp\left(-2\kappa k^2 t - \tfrac{2}{3}\kappa k^2 \alpha^2 t^3\right),
\]
\[
\iint_{\mathbb{T}^2} |\nabla \theta|^2\, d\bm{x} \leq (2\pi)^2(1+\alpha^2t^2)\exp\left(-2\kappa k^2 t - \tfrac{2}{3}\kappa k^2 \alpha^2 t^3\right).
\]
In contrast to Taylor’s dispersion scenario, where the enhancement scales inversely with $\kappa$, here we see that the enhanced rate of dispersion decreases as $\kappa$ decreases, making it especially significant in the limit of small diffusivity. 

\vspace{0.5cm}

\subsection{Enhanced dissipation through quantitative analysis} \label{s:resolvent}

In this section, we review recent mathematical results that quantify how shearing velocities influence the dynamics of \eqref{ad-eq}, primarily focusing on the infinite-channel setting originally studied by Taylor \cite{Taylor53}. We will present a selection of significant recent advances and outline the main methods employed to establish these results. 

We consider the scalar advection diffusion equation in the two-dimensional channel $\Omega=\mathbb{R}\times [0,1]$ with horizontal shear $(u(x_2),0)$. Then  the scalar function $\theta:\Omega\times(0,\infty)\rightarrow \mathbb{R}$ satisfies 
\begin{eqnarray}
    \partial_t \theta+u(x_2)\partial_{x_1}\theta&=&\kappa \Delta \theta \\
    \theta(x_1,x_2,0)&=&\theta_0(x_1,x_2)
\end{eqnarray}
where $u=u(x_2): D\rightarrow \mathbb{R}$ is a smooth function, representing the {\em shear profile} and $\kappa>0$ is again the molecular diffusivity. Here $\Delta=\partial_{x_1}^2+\partial_{x_2}^2$ is the classical Laplace operator. 
At the boundary of $\Omega=\mathbb{R}\times [0,1]$ we prescribe homogeneous Neumann boundary conditions 
\begin{equation}
\partial_{x_2}\theta(x_1, 0, t)=\partial_{x_2}\theta(x_1, 1, t) \mbox{ for all } x\in \mathbb{R } \mbox{ and } t>0. 
\end{equation}
We will see that most of the results on enhanced dissipation apply equally to the periodic channel $\mathbb{T}\times [0,1]$, the torus $\Omega=\mathbb{T}^2$, and can be  generalized to three dimensions.

Because of translation invariance in the horizontal variable, it is convenient to apply the Fourier transform in the $x_1$-variable, obtaining 
\begin{eqnarray}\label{eq-hat-theta}
    \partial_t \hat{\theta}+ik_1u(x_2)\hat{\theta}&=&\kappa (-k_1^2+\partial_{x_2}^2) \hat{\theta}, \\
    \hat{\theta}(k_1,x_2,0)&=&\hat{\theta}_0(k_1,x_2),
\end{eqnarray}
where $k_1\in \mathbb{R}$ is the horizontal wave-number.
The transformation
$$\hat{\theta}(k_1,x_2,t)=e^{-\kappa k_1^2 t}\zeta(k_1, x_2,t) $$
removes the term  $-\kappa k_1^2\hat{\theta}$ corresponding to horizontal diffusion, leading to the \textit{hypoelliptic} equation for the new function $\zeta$:
\begin{eqnarray}\label{new-eq}
    \partial_t \zeta+ik_1u(x_2)\zeta&=&\kappa \partial_{x_2}^2 \zeta \\
    \zeta(k_1,x_2,0)&=&\zeta_0(k_1,x_2)\,.
\end{eqnarray}
%We next consider the purely diffusive case $u=0$ in \eqref{new-eq}. 
due to the skew symmetry of ther advection term, testing the equation with $\zeta$ and integrating by parts gives
\begin{equation}\label{energy-es}
    \frac{d}{dt}\|\zeta(k_1, t)\|_{L^2}^2=-\kappa \|\partial_{x_2} \zeta(k_1, t)\|_{L^2}^2\,.
\end{equation}
Due to the Neumann boundary conditions on $\partial \Omega$, Poincaré's inequality cannot be directly applied in the $x_2$-direction. Nevertheless, we still obtain the following basic estimates:
\begin{equation}\label{decay-he}
   \|\zeta(k_1, t)\|_{L^2}^2 \leq \|\zeta_0(k_1)\|_{L^2}^2 
   \quad \Rightarrow \quad 
   \|\hat{\theta}(k_1, x_2)\|_{L^2}^2 \leq e^{-\kappa k_1^2 t}\|\hat{\theta}_0(k_1)\|_{L^2}^2\,.
\end{equation}
%where, we have reverted back to the function $\hat{\theta}$. 
This inequality shows that, for each fixed wavenumber $k_1$, the $L^2$-norm of $\hat{\theta}$ in the $x_2$-direction decays at most at a rate similar to that of the standard heat equation in the longitudinal direction, a bound that it is not optimal as it does not account for the influence of the shear flow’s advection on the dissipation.
Indeed, for each $k_1$, the dissipation timescale implied by the above estimate is comparable to that of the heat equation, on the order of $1/\kappa$.  More precisely, when $t \sim k_1^2/\kappa$, we loose the exponential decay, indicating that $\theta$ has mostly dissipated on that time scale.
%This can be seen from the exponential factor: when $t \sim k_1^2/\kappa$, the exponential essentially stabilizes, indicating that, at that time scale, the majority of the concentration is indeed dissipated.

In what follows, we will review recent results on enhanced dissipation in this context.
Over the past few years, two approaches have emerged as particularly effective in describing the enhanced dissipation effects due to the velocity field (but see \cite{BCZM24} for a different approach in the case of Hamiltonian flows). The first is the \emph{hypocoercivity method} introduced by Villani \cite{Villani09}, which systematically constructs modified energy functionals to capture subtle dissipative mechanisms. The second is a \emph{functional analytic method} based on resolvent estimates \cite{Wei18}. Both methods have provided deep insights into the interplay between advection and diffusion, yielding quantitative estimates on enhanced dissipation rates.
Under certain assumptions, we will show sharper estimates than those obtained by energy estimates as in \eqref{decay-he}, showing a time-decay of the form
\begin{equation}
    \|\hat{\theta}(k_1,t)\|_{L^2}^2 \leq e^{-\lambda_{\rm{ed}} t}\|\hat{\theta}_0(k_1)\|_{L^2}^2\,,
\end{equation}
where $\lambda_{\rm{ed}}$ depends only on $\kappa$ and $k_1$, and $\lambda_{\rm{ed}} \gg \kappa k_1^2$. This phenomenon illustrates how advection can substantially accelerate the decay of energy, encoding enhance dissipation, and  resulting in much faster convergence to equilibrium than one would expect from diffusion alone. In other words, the relevant timescale $t \sim 1/\lambda_{\rm{ed}}$ can be significantly shorter than the purely diffusive timescale $t \sim 1/\kappa$.

\subsubsection{Upper bounds via hypocoercivity}
In \cite{CZG23}, the authors assume that the velocity field $u$ is either monotone (no critical points) or belongs to $C^2$ with a finite number of non-degenerate critical points $\bar{x}_2^1, \ldots, \bar{x}_2^N$, where each satisfies
\begin{equation}\label{class1}
u''(\bar{x}_2^i) \neq 0 \quad \forall i=1,\ldots,N\,.
\end{equation}
This non-degeneracy ensures that the flow structure near each critical point is amenable to rigorous analysis, thereby allowing for precise quantification of enhanced dissipation rates.
To further limit the complexity, it is also assumed that the function $u$ has no critical points at the boundary $\partial{\Omega}$. Under these assumptions  Coti-Zelati and Gallay prove the folllowing result. 

\begin{theorem}[Theorem 1.1 \cite{CZG23}, simple critical points]\label{th-scp}
   Let $\Omega=\mathbb{R}\times[0,1]$ and $u\in C^m([0,1])$. Then there exists a positive constant $C$ such that for all $\kappa>0$, $k_1\neq 0$ and all initial data $\hat{\theta}_0$, the solution $\hat{\theta}$ to \eqref{eq-hat-theta} satisfies
    \begin{equation}\label{est-scp}
    \|\hat{\theta}(k_1, t)\|_{L^2(\Omega)}\lesssim e^{-\kappa k_1^2t-C\lambda_{\kappa, k_1}t}\|\hat{\theta}_0\|_{L^2(\Omega)},
    \end{equation}
    where 
    \begin{equation}
    \lambda_{\kappa, k_1}=\begin{cases}
                           \kappa^{\frac{1}{3}}|k_1|^{\frac{2}{3}} & \mbox{ if } 0< \kappa\leq |k_1|,\\
                            \frac{k_1^2}{\kappa}& \mbox{ if } 0< |k_1|\leq \kappa,
                          \end{cases}
    \end{equation}
if $u$ is monotone, i.e., $|u'|>0$, and 
  \begin{equation}
    \lambda_{\kappa, k_1}=\begin{cases}
                           \kappa^{\frac{1}{2}}|k_1|^{\frac{1}{2}} & \mbox{ if } 0< \kappa\leq |k_1|,\\
                            \frac{k_1^2}{\kappa}& \mbox{ if } 0< |k_1|\leq \kappa,
                          \end{cases}
    \end{equation}
if $u$ admits at most simple critical points, i.e., \eqref{class1} is satisfied.
\end{theorem}

In these results, one can observe a transition from an enhanced dissipation regime, occurring when 
$\kappa\leq |k_1|$, to a``strong diffusion" regime, which emerges when $ |k_1|\leq \kappa$. Because of the rate proportional to $\frac{1}{\kappa}$, the second regime is reminiscent of the one described by Taylor and for this reason one can name it ``Taylor dispersion regime". 

If we are interested in identifying the timescale at which the energy dissipation becomes effectively time-independent, we note that
%that the \emph{diffusive timescale} is $ t \sim \frac{1}{\kappa} $, while the \emph{enhanced dissipation timescale} is 
$ t \sim \frac{1}{\kappa^{1/3}} $ for monotone shear flows and $ t \sim \frac{1}{\kappa^{1/2}} $ for shear flows admitting simple critical points.

\begin{figure}
\centering
\includegraphics[scale=0.5]{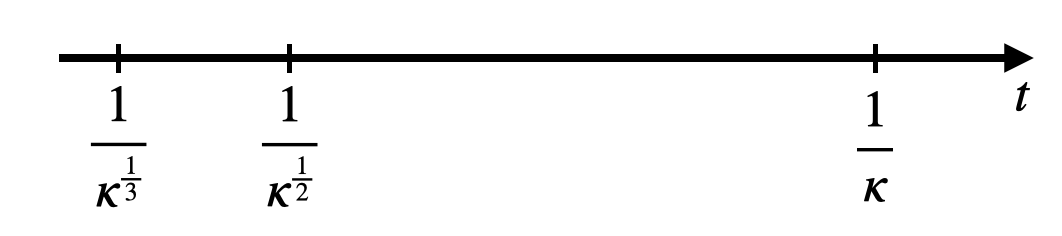}
\caption{Comparison between the pure diffusion time scale $t\sim \frac{1}{\kappa}$ and the enhanced dissipation time scales in the case of a monotone shear flow $t\sim\frac{1}{\kappa^{\frac{1}{3}}}$ and one admitting simple critical points  $t\sim\frac{1}{\kappa^{\frac{1}{2}}}$. }
\end{figure}
In the Taylor dispersion regime $|k_1|\ll \kappa$, we observe that the inequality
\[
\frac{k_1^2}{\kappa} \gg \kappa k_1^2
\]
holds when $\kappa \ll 1$, which in turn requires $|k_1| \ll 1$. Consequently, in \cite{CZG23}, the ``Taylor dispersion" regime is only observed for very small (real) wavenumbers (large wavelengths). A similar large-wavelength assumption was also made by G.~I.~Taylor (see Section~\ref{enhanced-dispersion}) when positing that the flow and scalar fields are sufficiently slowly varying, a key hypothesis ensuring that the dynamics can be reduced to \eqref{crucial-reduction}.

In \cite{BedrossianCotiZelati17} the authors considered $\Omega=\mathbb{T}$ and relaxed the hypothesis of the critical points, which can now be degenerate. They required that  the shear profile $u$ has a finite number of critical points, $\bar{x}_2^1,\cdots, \bar{x}_2^N$ such that for some $m_0\in \mathbb{Z}_+$,
\begin{equation}\label{class2}
\frac{d^{m_0}}{dx^{m_0}}(u')(\bar{x}_2^i)\neq 0\quad \forall i=1,\cdots, N\,.
\end{equation}
Thus $u$ is assumed of class at least $C^{m_0+1}$.
%where $m_0$ denotes the maximal order of derivative for which $u'$ is vanishing at the critical points. 
Strictly monotone shear flows correspond to the case $m_0=0$ in this definition.  
%Compared to the assumption \eqref{class1}, this class of functions includes points that are possibly \textit{degenerate}. 
On the other hand, the function $f(x)=x^4$ has a degenerate critical point at $x=0$, since $f'(0)=0$, and satisfies condition \eqref{class2} with  $m_0=3$.
%since $(f')'''(0)\neq 0$. %Notice that this type of critical points are not included in the class of functions for which \eqref{class1} hold.
Under these assumptions Bedrossian \& Coti-Zelati in \cite{BedrossianCotiZelati17} prove the following result.

\begin{theorem}[Theorem 1.3 \cite{BedrossianCotiZelati17}] \label{thm:enhancedShearChannel}
   Let $\Omega=\mathbb{T}\times[0,1]$ and $u\in C^{m_0+1}$ satisfies \eqref{class2}. Suppose that there exists $\eta_0\ll 1$ (depending only on $u$) such that if $k_1$ satisfies 
   $\kappa|k|^{-1}\leq \eta_0$ for some , then
   the solution $\hat{\theta}$ to \eqref{eq-hat-theta} satisfies
    \begin{equation}
    \|\hat{\theta}(k_1, t)\|_{L^2(\Omega)}\lesssim e^{-\kappa k_1^2t-C\lambda_{\kappa, k_1}t}\|\hat{\theta}_0\|_{L^2(\Omega)}
    \end{equation}
    where 
    \begin{equation}
    \lambda_{\kappa, k_1}=\frac{\kappa^{\frac{m_c+1}{m_c+3}}|k|^{\frac{2}{m_c+3}}}{(1+\log|k_1|+\log \kappa^{-1})^{2}}.
    \end{equation}
    where $m_c=\max\{m_0,1\}$.
\end{theorem}

In the same paper the authors prove a similar result on the two-dimensional torus.

\begin{theorem}[Theorem 1.1 \cite{BedrossianCotiZelati17}] \label{thm:enhancedShearTorus}
   Let $\Omega=\mathbb{T}^2$ and let $u\in C^{m_0+1}$ satisfy \eqref{class2}. Suppose that there exists $\delta_0\ll 1$(depending only on $u$) such that if  
   $\kappa|k|^{-1}\leq \delta_0$ for, then
   the solution $\hat{\theta}$ to \eqref{eq-hat-theta} satisfies
    \begin{equation}\label{est-torus}
    \|\hat{\theta}(k_1, t)\|_{L^2(\Omega)}\lesssim e^{-\kappa k_1^2t-C\lambda_{\kappa, k_1}t}\|\hat{\theta}_0\|_{L^2(\Omega)}
    \end{equation}
    where 
    \begin{equation}
    \lambda_{\kappa, k_1}=\frac{\kappa^{\frac{m_0+1}{m_0+3}}|k_1|^{\frac{2}{m_0+3}}}{(1+\log|k_1|+\log \kappa^{-1})^{2}}.
    \end{equation}
\end{theorem}

The only difference between Theorems \ref{thm:enhancedShearChannel} and \ref{thm:enhancedShearTorus} is  the parameter $m_c=\max\{m_0,1\}$ instead of $m_0$ in the decay rate for the periodic channel. This difference is due to the presence of a boundary component, which  has an effect comparable to that of  an additional critical point  vanishing to first order.
We note that the logarithmic correction in $\lambda_{k_1,\kappa}$
  has a purely technical origin. In \cite{wei2019enhanced}, Wei and Zhang eliminated this logarithmic correction in their estimates by employing time-dependent weights. 
  %within the augmented energy functional.

These results were established using the \textit{hypocoercivity method}, a framework developed by Cédric Villani in his seminal work \cite{Villani09}. 
Villani's hypocoercivity framework has been successfully applied to a wide range of settings, particularly in kinetic theory, including the Boltzmann and Fokker–Planck equations, \cites{desvillettes2001trend, dolbeault2015hypocoercivity, herau2007short}. Hypocoercivity has proved an effective method to establish enhanced dissipation. Besides the results already recalled, we mention the articles \cites{CH24,CZDG-V23,H22,ABN22}
 To provide a clearer presentation, we will illustrate the method (or rather, the ``scheme") while omitting the hat-notation for simplicity.

\smallskip

\noindent{\em Sketch of proof. } 
The hypocoercivity scheme involves constructing an augmented energy functional, $\Phi$, that includes carefully chosen cross terms and weights. These additional components are designed to capture the subtle interplay between dissipative and non-dissipative parts of the system, thereby guaranteeing exponential convergence to equilibrium under appropriate conditions.

\smallskip
\noindent\textbf{Step 1: Defining the augmented energy functional.}  
The hypocoercivity scheme begins by introducing an augmented energy functional:
\[
\Phi = \frac{1}{2} \biggl[ \|\theta\|_{L^2}^2 + \alpha\|\partial_{x_2}\theta\|_{L^2}^2 
+ 2\beta \operatorname{Re}\langle i k_1 \theta\, \partial_{x_2}u, \partial_{x_2}\theta \rangle 
+ \gamma k_1^2 \|\theta \partial_{x_2}u\|_{L^2}^2 \biggr],
\]
where $\alpha,\beta,\gamma$ are parameters to be determined later.

We then compute the time derivative $\frac{d}{dt}\Phi$. At this stage, boundary conditions are crucial. Nevertheless, for common settings such as the torus, or periodic or infinite channels with Neumann boundary conditions, all boundary terms arising from integration by parts vanish. After some straightforward manipulations, one obtains
\[
\frac{d}{dt}\Phi \leq -G + B,
\]
where $G > 0$ represents the ``good'' terms directly related to $\Phi$, and $B > 0$ corresponds to the ``bad'' error terms. The next steps will involve judiciously choosing $\alpha,\beta,\gamma$ to ensure that $G$ dominates $B$, thereby guaranteeing exponential decay of $\Phi$.

\smallskip

\textbf{Step 2: Partition of unity}. 
This step is unnecessary under the assumption of simple critical points. However, for general flows satisfying \eqref{class2}, it is required because the choice of $\alpha, \beta$, and $\gamma$ depends on the local behavior of $u$ near critical points, which may be degenerate.  

To handle this situation, we group all critical points, $\bar{x}_2^1, \cdots, \bar{x}_2^N$, of the same order: For each $j = 0, \dots, m_0$, we define  
\[
\phi_j(x_2) = \sum_{i \in E_j} \tilde{\phi}_i(x_2),
\]  
where $\tilde{\phi}_i$ is a smooth function equal to 1 in a small ball of radius $\delta$ around the $i$-th critical point and vanishes outside a ball of radius $2\delta$, and $\tilde{\phi}_0(x_2)=1-\sum_{i=1}^{N}\tilde{\phi}_i(x_2)$. The collection of $\{\tilde{\phi}_i\}_{i=0}^{N}$ forms a partition of unity. The set $E_j$ consists of indices $i$ such that  
\[
u^{(l)}(\bar{x}_2^i) = 0 \quad \forall l = 1, \dots, j, \quad u^{(j+1)}(\bar{x}_2^i) \neq 0.
\]

\smallskip

\textbf{Step 3: Spectral gap estimate}.
In the localized version, one uses a spectral gap estimate of the form  
\[
\sigma^{\frac{j}{j+1}} \|f_j\|_{L^2}^2 \lesssim \sigma \|\partial_{x_2} f_j\|_{L^2}^2 + \|u' f_j\|_{L^2}^2,
\]  
which is essential for controlling the error terms in the analysis.
\smallskip

\textbf{Step 4: Choice of balancing parameters and closing argument}.
In the case of simple critical points, there exist constants $\alpha = \alpha(\kappa, k_1)$, $\beta = \beta(\kappa, k_1)$, and $\gamma = \gamma(\kappa, k_1)$ such that  
\[
\frac{d}{dt} \Phi + C \nu^{\frac{1}{2}} |k|^{\frac{1}{2}} \Phi \leq 0.
\]  
Applying Grönwall's inequality and standard semigroup estimates yields the desired bound \eqref{est-scp}.  

Under the more general assumption on critical points \eqref{class2}, the strategy remains similar, with the key difference that $\alpha, \beta$, and $\gamma$ are no longer constants but depend non-trivially on $x_2$ due to the varying behavior at different critical points. Specifically,  
\[
\alpha(x_2) = \sum_{j=0}^{m_0} a_j \phi_j(x_2), \quad \beta(x_2) = \sum_{j=0}^{m_0} b_j \phi_j(x_2), \quad \gamma(x_2) = \sum_{j=0}^{m_0} c_j \phi_j(x_2).
\]  
%\vspace{1cm}
The regularity and structure of the velocity field $u(x_1,x_2) = (u(x_2),0)$ have played a crucial role in all the results discussed above. In the enhanced dissipation phenomena considered in various settings---the torus $\mathbb{T}^2$, periodic channel $\mathbb{T}\times [0,1]$, and infinite channel $\mathbb{R}\times [0,1]$---the velocity $u(x_2)$ is assumed to be a bounded, sufficiently regular function with only finitely many (possibly degenerate) critical points. As we have seen, the maximal order of degeneracy of these critical points (as described in \eqref{class2}) significantly influences the decay rate. Indeed, critical points are essential since their associated gradients drive the mixing process.
Notably, the enhanced dissipation rate does not depend on the amplitude of $u$ (i.e., the difference between its maximum and minimum values). Furthermore, we have seeen that in the case of infinite channels where the passive scalar is subject to homogeneous Neumann boundary conditions in the $x_2$-direction, the hypocoercivity method produces an estimate that recovers behavior analogous to Taylor’s shear dispersion as $k_1 \to 0$.

 If homogeneous Dirichlet boundary conditions were imposed instead, the Poincaré inequality would immediately yield a sharper upper bound compared to the ``Taylor dispersion" upper bound given in \eqref{est-scp}. Specifically, applying the Poincaré inequality in \eqref{energy-es}, we obtain  
\[
\frac{d}{dt}\|\zeta(k_1, t)\|_{L^2}^2 = -\kappa C_p \|\zeta(k_1, t)\|_{L^2}^2,
\]  
which leads to the exponential decay estimate:  
\[
\|\hat{\theta}(k_1, t)\|_{L^2}^2 \leq e^{-\kappa k_1^2 t - C_p \kappa t} \|\hat{\theta}_0\|_{L^2}^2.
\]  
Clearly, the decay rate here is faster than the rate given by \(\exp(-\kappa k_1^2 t - \frac{k_1^2}{\kappa} t)\) in the regime \( |k_1| \leq \kappa \), since one can always arrange that
%. Indeed, comparing with the rate in \eqref{est-scp}, we observe that, up to allowing the constant \(C\) to worsen slightly, the inequality  
\[
C \frac{k_1^2}{\kappa} \leq C_p \kappa
\]  
for $ |k_1| \leq \kappa $, by taking $C$ slightly smaller.

In \cite{CZ-Dolce2020}, the authors posed the natural question of whether enhanced dissipation phenomena can also be mathematically demonstrated in the whole space $\mathbb{R}^2$. This question is particularly intriguing because, in the whole space, there is no dissipation mechanism arising from boundaries, and the Poincaré inequality no longer applies. The result is established for a specifically chosen velocity field:
\begin{equation}
\bm{u}(x_1,x_2)=(x_1^2+x_2^2)^{q} \begin{pmatrix}
   -x_2\\ x_1
\end{pmatrix}
\end{equation}
with some fixed $q\geq 1$ and mean zero initial data $\theta_0$. Passing to polar coordinates $(r,\phi)$ the transport term transforms simply to
\begin{eqnarray*}
    (\bm{u}\cdot \nabla)\theta 
    %& =r^{q+1}\left[ -\sin\phi[\partial_r f \cos\phi-\frac{1}{r}\sin\phi \partial_{\phi} \theta]+\cos\phi(\partial_r f \sin\phi-\frac{1}{r}\cos\phi \partial_{\phi} \theta)\right] \\
    =r^{q}\partial_{\phi}\theta
\end{eqnarray*}
and the dynamics are effectively described by the one-dimensional equation
\begin{eqnarray}\label{unidim-shear}
    &&\partial_t \theta +r^q \partial_{\phi}\theta=\nu \Delta \theta \quad \mbox{ in } r\in[0,\infty), t\geq 0\,, \\
    && \theta|_{t=0}=\theta_0  \quad \mbox{ in } r\in[0,\infty)\notag\,.
\end{eqnarray}
Given the periodicity in the angular variable $\phi$, advection does not act on functions that are constant in $\phi$. Therefore, the mean in $\phi$ must be subtracted to see enhancement.  We therefore set
$$\langle \theta \rangle(r,t)=\frac{1}{2\pi}\int_{\mathbb{T}} \theta (r,\phi, t)\, d\phi.$$
The $\phi$-average of $\theta$ satisfies the heat equation in $\mathbb{R}^+$, which easily gives polynomial decay 
$$\|\langle\theta\rangle_{\phi}\|_{L^{\infty}}\lesssim \frac{C_0}{\sqrt{\kappa t}}\|\langle\theta_0\rangle_{\phi}\|_{L^2}\,.$$
Furthermore, since the advection term has an unbounded drift, suitable weights should be used. To this end, we  introduce the weighted $L^2$ norm:
$$\|f\|_{X}:=\int_0^{\infty}\int_{\mathbb{T}}(1+r^{2(q-1)})|g(r,\phi)|^2\,r\, dr\, d\phi\,.$$
Then Coti-Zelati \& Dolce obtained the following result.

\begin{theorem}[Theorem 1.1 \cite{CZ-Dolce2020}] For all $t\geq 0$
    $$\|\theta(t)-\langle\theta\rangle_{\phi}\|_{X}\lesssim e^{-\varepsilon_0\lambda_{\kappa}t}\|\theta_0(t)-\langle\theta_0\rangle_{\phi}\|_{X}\,,$$
    where $$ \lambda_{\kappa}=\frac{\kappa^{\frac{q}{q+2}}}{1+\frac{2(q-1)}{q+2}|\ln \kappa|}\,.$$
\end{theorem}

%Here the following compact notations is adopted: $$\langle \theta \rangle(r,t)=\frac{1}{2\pi}\int_{\mathbb{T}} \theta (r,\phi, t)\, d\phi$$
%is the average in the angular variable and 
%$$\|f\|_{X}:=\int_0^{\infty}\int_{\mathbb{T}}(1+r^{2(q-1)})|g(r,\phi)|^2\,r\, dr\, d\phi\,,$$
%is an appropriately weighted $L^2$ norm.

%Since the $\phi$-average of $\theta$ satisfies the heat equation in $\mathbb{R}^+$, one can easily show polynomial decay 
%$$\|\langle\theta\rangle_{\phi}\|_{L^{\infty}}\lesssim \frac{C_0}{\sqrt{\kappa t}}\|\langle\theta_0\rangle_{\phi}\|_{L^2}\,.$$

The resulting picture in this setting mirrors the findings of Theorem \ref{th-scp}: the velocity field $ u $ significantly enhances dissipation. This is evident from the enhanced dissipation time scales, $ t \sim \kappa^{-\frac{1}{3}} $ for $ p = 1 $ and $ t \sim \kappa^{-\frac{1}{2}} $ for $ p = 2 $, which are substantially shorter than the purely diffusive time scale, $ t \sim \kappa^{-1}, $ when $\kappa\ll 1$.
From this estimate we deduce the following picture. Initially, the velocity field $u$ enhances dissipation by mixing the scalar field, rapidly reducing its gradients along streamlines. As time progresses, and dissipation reaches the 
$\phi$-independent state, the scalar becomes nearly uniform along streamlines, leaving molecular diffusion to dissipate any remaining gradients perpendicular to the streamlines.
Diffusion ``takes over" because the enhanced dissipation mechanisms have largely completed their role, and only small-scale variations across streamlines remain to be smoothed out.

\subsubsection{Upper bounds via resolvent estimates}\label{Sec:ad-res-est}

In \cite{feng2023enhanced}, the authors of this review extended the approach in \cite{CZ-Dolce2020} to more general flows, namely,  circularly-symmetric flows with velocity fields of the form  
\begin{equation} \label{eq:R2Shears}
\bm{u}(x_1, x_2) = u\left(\sqrt{x_1^2 + x_2^2}\right) 
\begin{pmatrix}
    -x_2 \\
    x_1 \\
\end{pmatrix},
\end{equation}
where $u: [0, \infty) \to \mathbb{R}$ is a smooth profile with a finite number of (possibly degenerate) critical points, subject to certain conditions. Under these assumptions, the authors demonstrated enhanced dissipation at a rate $\lambda_\kappa = \kappa^{\frac{m}{m+2}}$, where $m$ quantifies the maximum order of degeneracy of the critical points.

In the periodic pipe $\Omega = \mathbb{D}(0, 1) \times \mathbb{T}$, where $\mathbb{D}(0, 1)$ denotes the unit disk in the $xy$-plane centered at the origin, the same authors considered a parallel pipe flow of the form  
\[
\bm{u}(x_1, x_2, x_3) = u(\sqrt{x_1^2+x_2^2})
\left(
    \sin(2\pi r)
    \begin{pmatrix}
        -x_2 \\
        x_1 \\
        0
    \end{pmatrix}
    +
    \cos(2\pi \sqrt{x_1^2+x_2^2})
    \begin{pmatrix}
        0 \\
        0 \\
        1
    \end{pmatrix}
\right),
\]
where $u: [0, 1] \to \mathbb{R}$ is again a smooth profile with a finite number of possibly degenerate critical points. Using cylindrical coordinates and assuming Neumann boundary conditions, the authors derived a quantitative estimate for the $L^2$-norm of the deviation $\theta - \langle \theta \rangle_{\phi, z}$, showing enhanced dissipation at the same rate, $\lambda_\kappa = \kappa^{\frac{m}{m+2}}$. Here, again $\langle \theta \rangle_{\phi, z}$ denotes the mean in the non-radial variables. Notably, these results were not obtained using the hypocoercivity framework but were achieved through resolvent estimates.

It is well-known that the norm of a semigroup can be bounded through the resolvent of the generator. 
The use of resolvent estimates offers the advantage of being a versatile method, applicable to fractional and higher-order operators. However, its applicability is restricted to autonomous flows that meet specific conditions, such as the flows studied by the authors.
We follow the approach in \cite{Wei18Mix} (see also \cites{MM22}, which utilizes a  Gearhart-Pr\"uss-Greiner type theorem  (we refer to \cite{EngelNagel00} for more details ), specialized to {\em maximally or $m$-accretive operators}. Informally, an unbounded operator $L$ on a Banach space $X$ with dense domain $\mathcal{D}(L)$ is accretive, if $-L$ is dissipative (see e.g. \cite{Pazy83}). In particular, the spectrum lies in the positive real line, so that the resolvent $ R(L;\lambda)= (I+\lambda \, L)^{-1}$, with $I$ the identity operator,  is a bounded operator on $X$ for all $\lambda>0$. The operator is then called maximally accretive if $R(L;\lambda)$ is surjective. 

For an $m$-accretive operator on a Hilbert space $H$, Theorem 1.3 in~\cite{Wei18Mix} gives that
\begin{align} \label{eq:WeiEstimate}
\|e^{-t L}\|_{\text{op}} \leq e^{-t\Psi( L)+\frac{\pi}{2}}\,,\quad \forall t\geq 0\,,
\end{align}
where $\|\cdot\|_{\text{op}}$ denotes the operator norm, and $\Psi( L)$ is the spectral function associated to the operator $L$ defined as:
\begin{equation} \label{}
 \Psi(L)=\inf\left\{\|( L-i\lambda)g\|: g\in \mathcal{D}(L),\,\lambda\in \mathbb{R},\, \|g\|=1\right\}\,.
\end{equation}
We apply this result when $L=L_\kappa:=\kappa \, (-\Delta)^\gamma - \bm{u}\cdot \nabla$, $\gamma=1,2$, and $\bm{u}$ is a shear flow on the 2D torus $\mathbb{T}^2$, a circularly symmetric flow in $\mathbb{R}^2$, or a 3D parallel pipe flow. 
Circularly symmetric flows in $\mathbb{R}^2$ can be interpreted as angular shears with radial profiles. Pipe flows, while not strictly unidirectional due to their helical streamlines, exhibit a cross-sectional component that remains circularly symmetric, while the vertical component is purely radial.
More general diffusion operators can also be considered, but for simplicity we only treat the case of  the bi-harmonic operator $\Delta^2$  and $\bm{u}=(u(x_2),0)$ a horizontal shear on $\mathbb{T}^2$ with {\em profile} $u$. The case of circularly symmetric flows and pipe flows can be treated in a similar fashion.
%For ease of notation, we denote a point $\bm{x}$ in $\mathbb{T}^2$ as $(x,y)$, as opposed to $(x_1,x_2)$.

Since the advection operator has a non-trivial kernel consisting of functions that are constant in the direction of the shear, where no enhanced dissipation can occur,  we work on the $L^2$-orthogonal complement of the kernel, which consists of functions with zero average in the direction of the shear. We denote this closed subspace of $L^2$ by $\mathring{L}^2$. We consider the operator $L$ as an unbounded operator on $\mathring{L}^2$.
By \eqref{eq:WeiEstimate}, one can estimate the rate of decay of the solution operator of the advection-hyperdiffusion equation $(\partial_t + L_\kappa)\theta = 0$ by estimating the spectral function. To this end, it is convenient to apply the Fourier Transform in the direction orthogonal to the shear and work with the transformed operator $L_{\kappa,k}:= \kappa \,(-\Delta_k)^2 + i k u(x_2)$, $\Delta_k:=\partial^2_{x_2} - k^2$, at {\em fixed wavenumber}  $k\ne 0$, which is an operator on functions of $x_2$ only.  
By Plancherel's theorem, it is enough to obtain spectral estimates for $L_{\kappa,k}$. Moreover, we have
\[
    \|e^{-t L_{\kappa,k} }\|_{\text{op}} \leq C\, \|e^{-t \tilde{L}_{\kappa,k} }\|_{\text{op}}
\]
where the constant $C$ is independent of $\kappa$ and $k$, and $\tilde{L}_{\kappa,k} := \kappa\, \partial^{2}_{x_2} +i k u(x_2)$. It can be shown \cite{CZDFM21} that  $L_{\kappa,k}$ is $m$-accretive on $L^2([0,2\pi])$ and hence \eqref{eq:WeiEstimate} applies. 
The resolvent estimate can be obtained via suitable energy estimates in the setting we consider and show enhanced dissipation provided a certain condition is imposed on the shear profile.

This condition was introduced in \cite{CZG23} to study both enhanced dissipation and Taylor dispersion in the higher-dimensional setting for the standard diffusion operator. It was then used by some of the authors of this work to obtain enhanced dissipation for the bi-harmonic operator, motivated by applications to the Kuramoto-Sivashinsky equation \cite{CZDFM21}. In Section \ref{s:applications} we discuss the implications of enhanced dissipation for the behavior of nonlinear systems.

\begin{assumption} \label{a:flow}
There exist $m, N \in \mathbb{N}$, $c_0>0$ and $\delta_0\in (0, 1)$, with the property that for any $\lambda\in \mathbb{R}$, $\delta\in (0, \delta_0)$ and $k\in \mathbb{Z}$ with $k\ne 0$, there exist $n\leq N$ and points $\{y^1,\cdots, y^n\}\in \mathbb{T}$ such that
\begin{align}\label{e:lowerbound}
| u(x_2)-\lambda| > c_0\,\delta^m\,,
\end{align}
for any $|x_2-y^j|\geq \delta$, $j\in\{1,\cdots,  n\}$.
\end{assumption}

An example of profile that satisfies Assumption  \ref{a:flow} above is that of the {\em Kolmogorov flow}, $u(x_2)=\sin(x_2)$, where $m=2$. More generally, we can take  $u(x_2)=\sin(x_2)^m$, $m\in \mathbb{N}$.

Under this assumption, one can estimate the spectral function associated to the operator $\tilde{L}_{\kappa,k}$, and consequently, obtain decay estimates for the semigroup as follows.

\begin{theorem}[ Proposition 3.1, Corollary 3.2 \cite{CZDFM21} ] \label{thm:ResolventEstimate}
 Let $\kappa\,|k|^{-1}\leq 1$, and let $P_k$ denote the $L^2$ projection onto the $k$-th Fourier horizontal mode.
There exists $\epsilon'_0>0$, independent of $\kappa$ and $k$, such that for every $t\geq 0$,
\begin{align} \label{eq:DecayRate}
\|e^{-L_\kappa t} \, P_k\|_{\text{op}} \leq e^{- \epsilon'_0 \kappa^\frac{m}{m+4}|k|^\frac{4}{m+4} t +\pi/2}.
\end{align}
In particular, $L_\kappa$ generates an exponentially stable semigroup in $\mathring{L}^2(\mathbb{T}^2)$ with rate:
\begin{equation} \label{eq:sembound}
   \|e^{-L_\kappa t}\|_{\text{op}} \leq  e^{-\lambda'_\kappa t +\pi/2},  \qquad t>0,
\end{equation}
where $\lambda'_\kappa=\epsilon'_0 \kappa^\frac{m}{m+4}$.
\end{theorem}

The term $\lambda'_\kappa$ characterizes the enhancement of dissipation in the sense that $\lambda'_\kappa / \kappa \to 0$ as $\kappa \to 0$. Specifically, on the torus, where $m \geq 0$ and any shear profile necessarily has at least one critical point, $\lambda'_\kappa = O(\kappa^{1/2})$.

Assumption \ref{a:flow} can be interpreted in an indirect way as a condition on the order of the critical points of  $u$. Indeed, a weaker estimate holds for any shear profile $u$ that has finitely many critical points of finite order, as such flows are mixing on the complement of the kernel of the advection operator (see e.g. \cite[Corollary 2.3]{CZDG20} and Section \ref{s:KSE}). 

\subsubsection{Upper bounds via Fourier splitting and Green functions}

We emphasize once again that, in the results discussed above for the whole space $\mathbb{R}^2$, the velocity fields were unbounded with respect to the spatial variables.
%(see Assumption 1). 
The growth at infinity of the velocity, which is assumed at least linear,  is a critical factor in achieving enhanced dissipation and appears to be independent of the specific method used to derive the decay estimates. Although we lack a formal proof, it seems unlikely that enhanced dissipation in $\mathbb{R}^n$ can be achieved with ``physical" velocity fields. In particular, if the profile $u$ in \eqref{eq:R2Shears} is required to be $L^p$-integrable for some $1 \leq p < \infty$, we are unable to demonstrate enhanced dissipation.

In \cite{Nobili-Pottel22}, the third author of this paper and Pottel proved quantitative bounds for mixing, quantified by the filamentation length $\Lambda (t)$ defined in \cite{MilesDoering18} as the ratio between the $H^{-1}$ norm and the $L^2$ norm of the scalar concentration $\theta$, in $\mathbb{R}^n$ with $n=2,3$.
Their estimates rely on the following upper and lower bound.

\begin{theorem}[Proposition 2.7 \& Theorem 1.1 \cite{Nobili-Pottel22}]
Let $n=2,3$, $\theta_0\in L^2(\mathbb{R}^n)$ be the initial condition,  $u(\cdot, t)\in L^2(\mathbb{R}^n)$  be divergence free and let  
$$r^{\ast}=r^{\ast}(\theta_0)=\sup\{r\in (-\frac{d}{2}, \infty)| \lim_{\delta \rightarrow 0} \delta^{-2r-d}\int_{|\xi|\leq \delta}|\hat{\theta}_0(\xi)|\, d\xi=0\}$$
denote the decay character of the initial data.
\begin{enumerate}
\item[i)]
 If  $r^{\ast}\in (-\frac{d}{2}, \infty)$ and 
\begin{equation}\label{id-assump}
\|u(t)\|_{L^2}\sim (1+t)^{\alpha}  \qquad \mbox{for some } \alpha>\frac 12-\frac d4\,. 
\end{equation}
Then there exists a positive constant $C$ depending on $d, r^{\ast}, \|\theta_0\|_{L^2}$ and $\alpha$ such that 
\begin{equation*}
    \|\theta(t)\|_{L^2}\lesssim \kappa^{-\max\{\frac d4+\frac{r^{\ast}}{2}, m\}}(1+t)^{-\min \{\frac d4+\frac{r^{\ast}}{2}, \frac{d}{4}+\frac 12\}}
\end{equation*}
for some $m\geq \frac{d}{2}+1$ depending on $\ceil{\frac{1}{|\alpha|}}$.
\item[ii)]
Under the more restrictive assumptions that  $r^{\ast}\in (-\frac{d}{2}, 1)$ and $$\|u(t)\|_{L^2}\sim (1+t)^{\alpha} \quad \mbox{ for some } \alpha>\frac{r^{\ast}}{2}+\frac 12\,,$$ if the time is sufficiently large, then 
\begin{equation}
    \|\theta(t)\|_{L^2}\gtrsim(\kappa(1+t))^{-\frac d4-\frac{r^{\ast}}{2}}\,.
\end{equation}
\end{enumerate}
\end{theorem}

How large time must be this is quantified explicitly in \cite{Nobili-Pottel22}.
As shown, the parameter that depends on the initial data, $r^{\ast}$, plays a pivotal role in determining the long-time decay of the solution’s energy. The proof of this result is based on the Fourier-splitting techniques introduced by Schonbek \cite{schonbek1986large}, which have been instrumental in establishing large-time asymptotics for the Navier–Stokes equations in unbounded domains.

These findings highlight a form of enhanced dissipation in the whole space that arises from the properties of the initial data. More specifically, the decay in the $L^2$ norm of the velocity field is governed by a parameter $\alpha$, which not only influences the well-posedness of the problem but is also crucial for closing the necessary estimates. In particular, the choice of initial data allows one to select a suitable stirring field that accelerates the decay beyond the standard polynomial rate $(\kappa t)^{-d/2}$ associated with pure heat diffusion. However, due to the perturbative nature of the Fourier-splitting method, quantifying the effect of “optimal” velocity fields remains challenging. As a result, the observed enhanced dissipation appears primarily dictated by the characteristics of the initial data.

In the diffusion-dominated regime described in (ii), the rate $(\kappa t)^{-d/4 - r^{\ast}/2}$ is sharp. It is worth noting that, in the whole space setting, the notion of separating time scales—so central to enhanced dissipation results in other contexts—does not seem to apply in the same straightforward manner.

One can ask how to  derive quantitative estimates that capture the influence of the velocity field on the whole space. As we have seen, the hypocoercivity method requires strong symmetry conditions: for example, the result in \cite{CZ-Dolce2020} applies only to a carefully chosen, radially symmetric velocity field that reduces the problem to that for a one-dimensional (weighted) shear. Similarly, resolvent estimates also depend on specific choices of the velocity field; in particular, we have seen that the results in \cite{CZ-Dolce2020} and \cite{feng2023enhanced} rely on the fact that the velocity field grows at infinity in the radial direction. In contrast, the Fourier splitting method, which can be employed on the entire space, accommodates a broad class of physical velocity fields but remains inherently perturbative.

A natural approach on $\mathbb{R}^n$ is to seek estimates through Green’s functions. Although universal Green’s functions for the advection–diffusion equation under minimal assumptions (that the velocity field is divergence-free and satisfies \eqref{id-assump}) are not available, upper and lower bounds have been established using methods for parabolic equations \cite{aronson1967bounds}. In \cite{carlen1996optimal}, for instance, assuming that the advecting field $\bm{u}$ is divergence-free and $\|\bm{u}(\cdot, t)\|_{L^{\infty}}\leq B(t)$, the authors derived an upper bound on the fundamental solution of the form:
\begin{equation}
G(\bm{x}, \bm{y}; T)\leq \frac{C_1}{T^{\frac{n}{2}}}\exp\left(-\frac{|\bm{x}-\bm{y}|^2}{T}\left(1-\frac{1}{|\bm{x}-\bm{y}|}\int_0^T B(\tau),d\tau\right)\right),
\end{equation} 
 for any $\bm{x},\bm{y} \in \mathbb{R}^n$. If, additionally, the flow is incompressible (i.e., $\bm{u}$ divergence-free) and satisfies $\|\bm{u}(\cdot, t)\|_{L^{\infty}}\sim \tfrac{1}{t^{\frac{1}{2}}}$, then Maekawa’s result \cite{maekawa2008lower} provides the following lower bound: 
\begin{equation} 
G(\bm{x}, \bm{y}; T)\geq \frac{C_3}{T^{\frac{n}{2}}}\exp\left(-C_4\frac{|\bm{x}-\bm{y}|^2}{T}\right), \quad T > s \geq s_0.
\end{equation}

These upper and lower pointwise estimates for the fundamental solution, known as \textit{Aronson estimates} \cite{aronson1967bounds}, are crucial for understanding the quantitative behavior of the advection–diffusion equation in the whole space. However, while a lower bound on the Green’s function is informative, it does not directly guarantee lower bounds on the solution itself. An intriguing open question remains: whether  there exists a physically meaningful velocity field (for instance, one belonging to certain $L^p$ spaces in both space and time) that induces enhanced dissipation in the whole space.

%\vspace{0.5cm}

\subsubsection{Enhanced dissipation by mixing vector fields}
In incompressible fluids, stirring leads to the generation of smaller scales and induces mixing, while diffusion acting as the damping mechanism, efficiently suppresses these small-scale structures. Understanding the interplay between mixing and diffusion, and in addition how they affect turbulence, is vital in the study of fluid dynamics~\cites{DoeringThiffeault06, CKRZ08, LinThiffeaultDoering11, Thiffeault12}. In~\cites{FengIyer19, Feng19,CZDG20}, the interaction between mixing and diffusion is quantitatively studied (see \cite{S23} for an approach to enhanced dissipation using optimal transport). In~\cites{FengIyer19, Feng19}, on a closed Riemannian manifold $M$, a vector field $\bm{u}$ is called strongly mixing with rate function $h$, if for all mean-free functions $f,  g\in  H^1(M)$, the corresponding flow map $\varphi_{s,t}$ satisfies
\begin{align}\label{e:mix}
\langle \varphi_{s,t}(f), g\rangle \leq h(t-s)\|f\|_{\dot H^1}\|g\|_{\dot H^1}\,,
\end{align}
where $\langle,\rangle$ the $L^2$ inner product,  $h: [0, +\infty)\to [0, +\infty) $ is a decreasing function vanishing at infinity and the flow map $\varphi_{s,t}$ is defined by
\begin{align}
\partial_t\varphi_{s,t}=\bm{u}(\varphi_{s,t},t)\,,\quad \varphi_{s,s}=\mathbf{Id}\,.
\end{align}
The above definition quantifies the filamentation in the scalar due to the flow. By duality, equation~\eqref{e:mix} yields
\begin{align}\label{e:mix1}
\|f\circ\varphi_{s,t}\|_{H^{-1}}\leq h(t-s)\|f\|_{\dot H^1}\,.
\end{align}
Let $0<\lambda_1<\lambda_2\leq \cdots$ be the eigenvalues of the Laplacian operator on $L^2$, and let $\{e_i\}$ be the Hilbert basis of mean-free functions in $L^2$.  We denote by $\mathbf{P}_N$  the projection onto $\{e_1, e_2,\cdots, e_N\}$. Since $\bm{u}$ is divergence free, $\|\varphi_{s,t}(f)\|_{L^2}$ is preserved. Then equation~\eqref{e:mix1} yields
\begin{align}
\|(I-\mathbf{P}_N)f\circ\varphi_{s,t} \|_{L^2}^2&=\|f\|_{L^2}^2-\|\mathbf{P}_Nf\circ\varphi_{s,t}\|_{L^2}^2\\\nonumber
&\geq \|f\|_{L^2}^2-\lambda_N\|f\circ\varphi_{s,t}\|_{ H^{-1}}^2\\\nonumber
&\geq \|f\|_{L^2}^2-\lambda_Nh^2(t-s)\|f\|_{\dot H^1}^2\,,
\end{align}
Since $h$ is decreasing as time increases, from the above estimate, one can see that the energy is transferred to  high frequencies. More precisely, equation~\eqref{e:mix1} gives the rate at which transport by $\bm{u}$ transfers energy from the low frequency part to the high frequency one in $L^2$ sense. 
If  the mixing rate function $h$ is known, the rate of enhanced dissipation can then be estimated, as exemplified in the following result.

\begin{theorem}[Theorem 4.1.1 \cite{Feng19}]\label{t:debymix}
Let $\bm{ u}$ be a vector field that is mixing with rate function $h$, then the solution of equation~\eqref{ad-eq} satisfies
\begin{align}
\|\theta(t)\|_{L^2}\leq e^{-\kappa H(\kappa)(t-s)}\|\theta(s)\|_{L^2}\,,
\end{align}
for any $0\leq s\leq t$, where $H(\kappa)$ solves
\begin{align}\label{e:rate}
2H(\kappa)\;h\left(\frac{1}{64\sqrt{\kappa \|\nabla \bm{u}\|_{L^\infty}H(\kappa)}}\right)=1\,.
\end{align}
\end{theorem}

The above theorem gives a quantitative estimate on the dissipation enhancement based on the mixing rate function. 
Indeed,  the fact that $H$ satisfies~\eqref{e:rate} immediately implies $H(\kappa)\to \infty$ as $\kappa \to 0$, given that  the mixing rate function $h$ vanishes at infinity.
In particular, if the mixing rate function is exponentially decreasing, we have that
\begin{align}
\|\theta(t)\|_{L^2}\leq e^{-C(t-s)/|\ln\kappa|^2}\|\theta(s)\|_{L^2}\,,
\end{align}
where $C$ is a universal constant.
This rate is consistent with the estimate obtained in~\cite{CZDG20}. We remark that, if $\bm{u}$ is Lipschitz continuous, then mixing at an exponential rate is optimal (we refer again to the review article \cite{CZCIMSurvey} and references therein for a more-in-depth discussion of this point).

For (say horizontal) shear flows ${\bm u}(\bm{x})=(u(x_2),0)$ on $\mathbb{T}^2$,  we restrict again our attention to  mean-free functions with zero horizontal average. Then the method of stationary phase gives that the flow generated by ${\bm {u}}$ is mixing with rate function $h(t)=Ct^{-1/2}$ (see equation (1.8) in~\cite{BedrossianCotiZelati17}). Applying Theorem~\ref{t:debymix}, one has
\begin{align*}
\|\theta(t)\|_{L^2}\leq Ce^{-c\kappa^{4/5}(t-s)}\|\theta(s)\|_{L^2}\,.
\end{align*}
While Theorem 1.1 in~\cite{BedrossianCotiZelati17} guarantees that the energy decays at a much faster rate,
\begin{align}
\|\theta(t)\|_{L^2}\leq Ce^{-\frac{c\kappa^{1/2}(t-s)}{|\ln \kappa|^2}}\|\theta(s)\|_{L^2}\,.
\end{align}

The proof of Theorem~\ref{t:debymix} is mainly based on the Fourier splitting method. Multiplying both sides in~\eqref{ad-eq} by $\theta$ and integrating in time gives that
\begin{align}\label{e:L2energy}
\frac{d}{dt}\|\theta\|_{L^2}^2=-2\kappa\|\nabla \theta\|_{L^2}^2\,.
\end{align}
Next, there are two cases. Given a fixed $\lambda>0$, if  $\|\nabla \theta(t)\|_{L^2}^2 \geq \lambda\|\theta(t)\|_{L^2}^2$ for all $t\geq s\geq 0$, we immediately have $\|\theta(t)\|_{L^2}\leq e^{-\lambda\kappa(t-s)}\|\theta(s)\|_{L^2}$. On the other hand, if at some time $t_0$,  $\|\nabla \theta(t_0)\|_{L^2}^2 <\lambda\|\theta(t_0)\|_{L^2}^2$, then it is enhanced dissipation that gives the decay of the $L^2$ norm. To be more precise,  there is a time $\tau>0$ such that  $\|\theta(t_0+\tau)\|_{L^2}$ becomes sufficiently small, which is the content of the following lemma.

\begin{lemma}[Lemma 4.1.3 \cite{Feng19}]\label{l:de}
Let  $\lambda_N$  be the largest eigenvalue of the Laplace operator satisfying $\lambda_N\leq H(\kappa)$, where $H(\cdot)$ is defined in~\eqref{e:rate}. If $\|\nabla\theta(t_0)\|_{L^2}^2<\lambda_N\|\theta(t_0)\|_{L^2}^2$, then we have
\begin{align}\label{e:decaytau}
\|\theta(t_0+\tau)\|_{L^2}^2\leq \exp\left(-\frac{\kappa H(\kappa)\tau}{8}\right)\|\theta(t_0)\|_{L^2}^2\,,
\end{align}
for \,
%\begin{align}
$\tau=2h^{-1}\left(\frac{1}{2\lambda_N}\right)$.
%\end{align}
\end{lemma}

%ith these observations, we balance the effects of both cases and obtain an optimal decay rate that holds over all time.

We briefly discuss the proof of Lemma~\ref{l:de}. The energy estimate implies
\begin{align}
\|\theta(t_0+\tau)\|_{L^2}^2 =\|\theta(t_0)\|_{L^2}^2-2\kappa\int_{t_0}^{t_0+\tau}\|\nabla\theta(r)\|_{L^2}^2\,dr\,.
\end{align}
The choice of $\lambda_N$ gives that $\frac{H(\kappa)}{2}\lambda_N<H(\kappa)$ when $\kappa\ll 1$. Hence, ~\eqref{e:decaytau} follows
provided
\begin{align}\label{e:tmp1}
\int_{t_0}^{t_0+\tau}\|\nabla\theta(r)\|_{L^2}^2\,dr\ge \frac{\lambda_N\tau}{8}\|\theta(t_0)\|_{L^2}^2\,.
\end{align}
We prove~\eqref{e:tmp1} by contradiction. Therefore, we  assume that the opposite
 inequality holds: 
\begin{align}\label{e:converse}
\int_{t_0}^{t_0+\tau}\|\nabla\theta(r)\|_{L^2}^2\,dr< \frac{\lambda_N\tau}{8}\|\theta(t_0)\|_{L^2}^2\,.
\end{align}
Let $\varphi$ solve the  transport equation by $\bm{u}$ with initial data $\varphi(t_0)=\theta(t_0)$:
\begin{align}\label{e:trans}
\partial_t\varphi+\bm{u}\cdot \nabla \varphi=0\,.
\end{align}
We observe that
\begin{align}\label{e:tmp2}
\int_{t_0}^{t_0+\tau}\|\nabla\theta(r)\|_{L^2}^2&\geq \lambda_N\int_{t_0+\frac{\tau}{2}}^{t_0+\tau}\|(I-P_N)\theta(r)\|_{L^2}^2\,dr \nonumber\\
&\ge \frac{\lambda_N}{2}\int_{t_0+\frac{\tau}{2}}^{t_0+\tau}\|(I-P_N)\varphi(r)\|_{L^2}^2\,dr\nonumber\\
&\quad-\lambda_N\int_{t_0+\frac{\tau}{2}}^{t_0+\tau}\|(I-P_N)\left(\theta(r)-\varphi(r)\right)\|_{L^2}^2\,dr\nonumber\\
&\geq \frac{\lambda_N\tau}{4}\|\theta(t_0)\|_{L^2}^2-\frac{\lambda_N}{2}\int_{t_0+\frac{\tau}{2}}^{t_0+\tau}\|P_N\varphi(r)\|_{L^2}^2\,dr\\
&\quad-\lambda_N\int_{t_0}^{t_0+\tau}\|\theta(r)-\varphi(r)\|_{L^2}^2\,dr\nonumber
\end{align}
The fact that $\bm u$ is mixing with rate $h$ combined with the condition in Lemma~\ref{l:de} gives
\begin{align}\label{e:transmix}
\int_{t_0+\frac{\tau}{2}}^{t_0+\tau}\|P_N\varphi(r)\|_{L^2}^2\,dr&\leq \lambda_N \int_{t_0+\frac{\tau}{2}}^{t_0+\tau}\|P_N\varphi(r)\|_{H^{-1}}^2\,dr\leq \frac{\tau\lambda_N}{2} h(\frac{\tau}{2})^2\|\theta(t_0)\|_{H^1}^2\nonumber\\
&\leq \frac{\tau\lambda_N^2}{2} h(\frac{\tau}{2})^2\|\theta(t_0)\|_{L^2}^2\,.
\end{align}
We also use that $\theta$ and $\varphi$ are close in $L^2$ over the time interval $[t_0,t_0+\tau)$ (Lemma 4.1.4 in~\cite{Feng19}):
\begin{align}\label{e:difference}
&\int_{t_0}^{t_0+\tau}\|\theta(r)-\varphi(r)\|_{L^2}^2\,dr\nonumber\\
&\quad\leq 2\sqrt{2\kappa}\tau^{3/2}\|\theta(t_0)\|_{L^2}\left(2\|\nabla\bm u\|_{L^\infty}\int_{t_0}^{t_0+\tau}\|\theta(t)\|_{H^1}^2\,dt+\|\theta(t_0)\|_{H^1}^2\right)^{1/2}\nonumber\\
&\quad\leq 2\sqrt{2\kappa}\tau^{3/2}\|\theta(t_0)\|_{L^2}\left(\frac{\|\nabla\bm u\|_{L^\infty}\lambda_N\tau}{4}+\lambda_N\right)^{1/2}
\end{align}
By inserting~\eqref{e:transmix} and~\eqref{e:difference} into~\eqref{e:tmp2} and using assumption~\eqref{e:converse}, we have that
\begin{align}
\frac{\lambda_N \tau\|\theta(t_0)\|_{L^2}^2}{8}&\geq \frac{\lambda_N\tau\|\theta(t_0)\|_{L^2}^2}{4}-\frac{\tau}{4}\lambda_N^3h(\frac{\tau}{2})^2\|\theta(t_0)\|_{L^2}^2\\
&\quad-2\lambda_N\sqrt{2\kappa}\tau^{3/2}\|\theta(t_0)\|_{L^2}\left(\frac{\|\nabla\bm u\|_{L^\infty}\lambda_N\tau}{4}+\lambda_N\right)^{1/2}\,.
\end{align}
Finally, substituting the chosen values for $\tau$ and $\lambda_N$, yields a contradiction and  the lemma is proven.

% we relate the advection diffusion system~\eqref{ad-eq} with the corresponding transport system 
% \begin{align}\label{e:trans}
% \partial_t\phi+\bm{u}\cdot \nabla \phi=0\,.
% \end{align}
% Recall that the mixing property of the vector field $\bm{u}$ would transfer energy to the high frequency part and that the systems~\eqref{ad-eq} and~\eqref{e:trans} are close within a short period of time. In the advection diffusion equation, the diffusive operator then aids in the decay of high-frequency components and contributes to the overall energy dissipation.
% As a result, after a certain amount of time, a favorable decay rate can still be achieved. By balancing the effects of both cases, we can obtain an optimal decay rate that holds over all time.

Constructing exponentially mixing flows is a challenging task, and has been studied extensively in the dynamical systems literature ~\cites{Anosov67,Pollicott85,Dolgopyat98,Liverani04,ButterleyWar20, TsujiiZhang23}. Unfortunately, the question of whether smooth, time-independent or time-periodic, mixing flows exist on the torus is still an open problem, though advances have been made in recent years  (see for example~\cites{YaoZlatos17, AlbertiCrippaEA19,EZ19,MHSW22,ELM23} for deterministic constructions and ~\cites{BBPS22,BCZG22} in the more flexible and powerful random setting).

\section{Applications to non-linear systems } \label{s:applications}

Enhanced dissipation due to the combined effect of mixing by a background flow and diffusion can have a stabilizing effect in non-linear systems, a property that has been exploited to quench reactions, prevent blow-up, and suppress the growth of solutions for instance (we mention the works \cites{FannjiangKiselevEA06, HouLei09, BerestyckiKiselevEA10,KiselevXu16, BH18,FengFengEA20,IyerXuEA21}. A  well-studied model in this context is the Keller-Segel chemotaxis model. It is well known that the classical Keller-Segel system can experience blow-up in finite time when the spatial dimension is greater than one. Kiselev and Xu~\cite{KiselevXu16} considered the dissipation-enhancing flow introduced in~\cite{CKRZ08} and demonstrated that the solution of the advective Keller-Segel equation does not blow up in finite time, provided the amplitude of the  flow is large enough. Later in~\cite{IyerXuEA21}, flows with small dissipation times were studied and applied to the advective Keller-Segel system to ensure global well-posedness. For the generalized Keller-Segel equation with a fractional Laplacian and advection by a dissipation-enhancing flow, global well-posedness was discussed in~\cites{hopf2018aggregation,shi2020suppression}. In addition, the effect of advection by shear flows~\cite{BH18,He18} and planar helical flows~\cite{FengShiEA22} have also been studied in chemotaxis models, leading to suppression of  blow-up. In~\cite{zeng2021suppression}, the coupled Keller-Segel-Navier-Stokes system near the Couette flow $(x_2,0)$ was considered, establishingl well-posedness guaranteed if the amplitude of the flow is large enough. 

In this section, we mainly focus on two examples, suppression of phase separation in the Cahn-Hilliard equation ~\cite{FengFengEA20} and persistence of solutions for the 2D Kuramoto-Sivanshinsky equation  ~\cites{CZDFM21, FM22}.

\subsection{Mixing flows: applications to the Cahn-Hilliard equation} \label{s:CHE}
The  Cahn- Hilliard equation  \cites{CahnHilliard58, Cahn61} is a classic model of phase separation in binary mixtures, where the system spontaneously separates into two regions corresponding to different concentrations of the two components.  The evolution of the normalized concentration difference $c$ between the two phases is given by 
\begin{align}~\label{e:ch}
\partial_t c  +\gamma \Delta^2 c = \Delta(c^3-c).
\end{align}
Here $\sqrt \gamma$ is the Cahn number, which is related to the surface tension at the interface between phases. When $\gamma$ is small, solutions spontaneously form islands where $c=\pm 1$ separated by thin transition regions, a phenomenon that has been well studied, for instance in~\cites{ElliottSongmu86,Elliott89,Pego}.
In~\cite{FengFengEA20}, the effect of stirring on the spontaneous phase separation was investigated, by studing the following advective Cahn-Hilliard equation:
\begin{align}\label{e:ach}
\partial_t c +\bm{u}\cdot \nabla c +\gamma\Delta^2 c = \Delta(c^3-c)\,,
\end{align}
where $\bm{u}$ is some divergence-free vector field with the mixing property. For simplicity, we consider the above system on $\mathbb{T}^d$, with $d=2$ or $3$.
The goal is to show that if the stirring velocity field is sufficiently mixing, then there is no phase separation. More precisely, if the dissipation time of $\bm{u}$ is
small enough, then we can show that $c$ converges exponentially to the total concentration $\bar c = \int_{\mathbb{T}^d} c_0 dx$, where $c_0$ denotes the initial data. 

To present this result, we first define the concept of \textit{dissipation time} associated with the advection-hyperdiffusion equation,
\begin{align}\label{e:ahd}
\partial_t \theta +\mathbf u\cdot\nabla \theta+\gamma \Delta^2\theta=0\,.
\end{align}
The dissipation time, denoted by $t_{dis}(\bm{u},\kappa)$, is defined as the smallest time $t\geq 0$ such that the solutions to the above system~\eqref{e:ahd} satisfies
\begin{align} \label{eq:dissipationTime}
\|\theta(s+t)\|_{L^2}\leq \frac{1}{2}\|\theta(s)\|_{L^2}\,,
\end{align}
for all $s\geq 0$ and all mean-free initial data $\theta_0\in L^2$. A similar definition holds for the advection-diffusion equation~\eqref{ad-eq}. Then suppression of blow-up is a consequence of the following result.

\begin{theorem}[Theorem 1.2 \cite{FengFengEA20}]\label{t:ch}
Let $d=\{2,3\}$, $\bm{u}\in L^\infty((0,\infty);W^{1,\infty}(\mathbb{T}^d))$, and $c$ be the solution of~\eqref{e:ach} with initial data $c_0\in H^2(\mathbb{T}^d)$.\\
(1) When $d=2$, for any $\mu>0$, there exists a time $T_0=T_0(\|c_0-\bar c\|_{L^2}, \bar c, \gamma, \mu)$ such that if $t_{dis}(\bm{u},\gamma)< T_0$, then for every $t\ge 0$, we have
\begin{align}\label{e:cdecay}
\|c(t)-\bar c\|_{L^2}\leq 2 e^{-\mu t}\|c_0-\bar c\|_{L^2}\,.
\end{align}
 (2) When $d=3$, for any $\mu>0$, there exists a time $T_1=T_1(\|c_0-\bar c\|_{L^2}, \bar c, \gamma, \mu)$ such that if 
 \begin{align}\label{e:disass}
(1+\|\nabla \bm u\|_{L^\infty})^{1/2}t_{dis}(\bm{u},\gamma)<T_1\,,
 \end{align}
then~\eqref{e:cdecay} still holds for every $t\ge 0$.
\end{theorem}

To apply Theorem~\ref{t:ch}, one needs to produce vector fields $\bm{u}$ satysfying$t_{dis}(\bm{u},\gamma)< T_0$ when $d=2$, and $(1+\|\nabla \bm u\|_{L^\infty})^{1/2}t_{dis}(\bm{u},\gamma)<T_1$ when $d=3$. We accomplish this task by employing mixing flows with large amplitude and a sufficiently fast mixing rate $h$, a property that can be proved via a suitable rescaling of the flow. This approach has been studied  for example in ~\cites{CKRZ08, KiselevShterenbergEA08,CZDG20,FengIyer19,Feng19} for the standard Laplacian. With minor adjustments, the proofs  can be adapted to our context for the bi-harmonic operator.

\begin{proposition}[Proposition 1.4 \cite{FengFengEA20}]
Let $v\in L^\infty([0, \infty);C^2(\mathbb{T}^d))$, and define $\bm{u}_A(\bm{x},t)=Av(\bm{x},At)$. If $v$ is mixing with rate function $h$, then 
\begin{align*}
t_{dis}(\bm{u}_A,\gamma)\to 0\,,\text{~as $A\to \infty$}.
\end{align*}
Moreover, if the rate function $h$ satisfies
\begin{align*}
th(t)\to 0\,\text{~as $t\to \infty$},
\end{align*}
then
\begin{align*}
(1+\|\nabla \bm u\|_{L^\infty})^{1/2}t_{dis}(\bm{u}_A,\gamma)\to 0\,,\text{~as $A\to \infty$}.
\end{align*}
\end{proposition}

We remark here that the disadvantage of constructing flows with small dissipation times in this manner is that known examples of strongly mixing flows typically tend to be either quite complicated or lack sufficient regularity (see for instance~\cites{CKRZ08,YaoZlatos17,EZ19,AlbertiCrippaMazzucato19}). However, there are many flows that are not mixing but still have small dissipation times. Such flows are sufficient for our purposes as well.
For example, in~\cite{IyerXuEA21} flows with arbitrarily small dissipation times are constructed by rescaling a general class of smooth (time-independent) cellular flows. These flows have been employed in aggregation models to prevent the blow-up of solutions.

We now briefly outline the  proof of Theorem~\ref{t:ch}. In the case $d=2$, the result follows immediately from the two lemmata below. 

\begin{lemma}[Lemma 2.2 \cite{FengFengEA20}]\label{l:ch1}
For any $t_0\geq 0$, we have 
\begin{align}
\sup_{0\leq \tau\le \gamma\ln 2}\|c(t_0+\tau)-\bar c\|_{L^2}^2\leq 2\|c(t_0)-\bar c\|_{L^2}^2\,.
\end{align}
Moreover, if for some $\tau\in (0, \gamma\ln 2)$ and $\mu>0$ we have 
\begin{align}
\frac{1}{\tau}\int_{t_0}^{t_0+\tau}\|\Delta c\|_{L^2}^2\,ds \geq \frac{2+2\gamma\mu}{\gamma^2}\|c(t_0)-\bar c\|_{L^2}^2\,,
\end{align}
then
\begin{align}\label{e:largedecay}
\|c(t_0+\tau)-\bar c\|_{L^2}\leq e^{-\mu \tau}\|c(t_0)-\bar c\|_{L^2}\,.
\end{align}
\end{lemma}

This lemma establishes two key facts:
\begin{enumerate}
\item The quantity $\|c-\bar c\|_{L^2}^2$ can be controlled  within a certain time frame; 

\item If on this time frame, the average in time of $\|\Delta c\|_{L^2}^2$ is large enough,  the $L^2$ norm decays exponentially.
\end{enumerate}

When instead  the time average of $\|\Delta c\|_{L^2}^2$ is small, the next lemma shows that if $t_{dis}(\bm{u},\gamma)$ is small enough,  the variance of $c$ continues to decrease  by a constant fraction after time $t_{dis}(\bm{u},\gamma)$.

\begin{lemma}[Lemma 2.3 \cite{FengFengEA20}]\label{l:ch2}
For any $t_0\geq 0$, there exists a time
\begin{align}
T_0'=T_0'(\|c(t_0)-\bar c\|_{L^2},\bar c,\gamma,\mu)\in (0, \gamma\ln 2]\,,
\end{align}
such that if 
\begin{align}
t_{dis}(\bm{u}, \gamma)&\leq T_0'\,,\\
\frac{1}{t_{dis}(\bm{u}, \gamma)}\int_{t_0}^{t_0+t_{dis}(\bm{u}, \gamma)}\|\Delta c\|_{L^2}^2\,ds&\leq \frac{2+2\gamma \mu}{\gamma^2}\|c(t_0)-\bar c\|_{L^2}^2\,,
\end{align}
then~\eqref{e:largedecay} still holds at time $\tau=t_{dis}(\bm{u},\gamma)$. Moreover, the time $T_0'$ can be chosen to be decreasing as a function of $\|c(t_0)-\bar c\|_{L^2}$.
\end{lemma}

With these two lemmata at hand, Theorem~\ref{t:ch} in 2D  follows in a straightforward fashion. Lemma~\ref{l:ch1} can be proved by   energy estimates combined with Sobolev inequalities. For the proof of  Lemma~\ref{l:ch2},  we define $\mathcal{S}_{s,t}$ to be the solution operator of the system~\eqref{e:ahd}. For simplicity, we assume $t_0=0$. Then using Duhamel's principle, we have that
\begin{align}
c(t_{dis})-\bar c=\mathcal{S}_{0,t_{dis}}(c_0-\bar c)+\int_{0}^{t_{dis}}\mathcal{S}_{s,t_{dis}}(\Delta (c^3(s)-c(s)))\,ds\,.
\end{align}
By the definition of $t_{dis}$, the solution operator $\mathcal{S}_{s,t}$ halves the variance of $c$ in time $t_{dis}$. Furthermore, since $\mathcal{S}_{s,t}$ is an $L^2$ contraction, we obtain
\begin{align}
\|c(t_{dis})-\bar c\|_{L^2}\leq \frac{1}{2}\|c_0-\bar c\|_{L^2}+\int_{0}^{t_{dis}}\|\Delta (c^3(s)-c(s))\|_{L^2}\,ds
\end{align}
If $t_{dis}$ is small enough, one can finally prove that the nonlinear terms can be neglected in estimating the decay of the variance.

In the case $d=3$, the result follows in a similar fashion. However, to establish the analog of of Lemma~\ref{l:ch2}, a stronger assumption on 
$t_{dis}$ is necessary, specifically condition~\eqref{e:disass}.

\subsection{Shear flows: applications to the Kuramoto-Sivanshinsky equation} \label{s:KSE}

In this last section, we discuss applications of enhanced dissipation for shear and shear-like flows to the control of non-linear PDEs.
To complement the results already discussed, especially in Section \ref{s:CHE}, we will exemplify the effect of enhanced dissipation on the {\em Kuramoto-Sivashinsky equation} (KSE for short), a fourth-order nonlinear parabolic equation that models flame front propagation \cites{Kuramoto78,Sivashinsky77,Sivashinsky80}. The KSE is a paradigm for long-wave instability in dissipative system and, even in one space dimension, it exhibits spatio-temporal chaos.

In what follows, if $f$ is a function of space and time we will use the notation $f(t)$ to denote the function of $x$, $f(t)(\bm{x})=f({\bm{x}}, t)$.

We will discuss the KSE with periodic boundary conditions, a setting that can be justified in the context of wave propagation. For ease of notation, we continue to denote the torus with periods $L_i$, $i=1,\ldots, d$, by $\mathbb{T}^d$. Since it models the motion of an interface (the flame front), the physically relevant dimensions are $d=1,2$. 
Then the KSE in its primitive form is the following equation for a scalar potential $\varphi:[0,T)\times \mathbb{T}^d$, $0<T\leq \infty$:
\begin{equation} \label{eq:KSEscalar}
  \partial_t \varphi + (\Delta^2+\Delta) \varphi +|\nabla \varphi|^2=0,
\end{equation}
complemented by an initial condition $\\varphi(0)=\varphi_0$. It also has a vectorial form for ${\bm v}=\nabla \varphi$:
\begin{equation} \label{eq:KSEvector}
  \partial_t {\bm v} + (\Delta^2+\Delta) {\bm v} + {\bm v} \cdot \nabla {\bm v}=0, \qquad \text{curl\,} {\bm v}={\bm 0}.
\end{equation}
The linearized operator $\Delta^2+\Delta$ has finitely-many growing modes as soon as any period is strictly greater than $2\pi$, encoding the large scale instability. The non-linearity couples large and small scales where (hyper)dissipation is dominant, leading to the complex spatio-temporal dynamics observed in the system (for careful numerical simulations in $1$ and $2$D, we mention \cite{KSE1DNum,KSE2DNum} and references therein).

While the analysis of the 1D KSE is well-developed, due to the special structure of the non-linearity, which is an exact derivative, although important questions remain open on the characteristic long-time behavior of solutions, especially with respect to the dependence on the size of the periodic box, much less is known in the two-dimensional case. 
%As a matter of fact, the global-in-time existence of solutions is still essentially open, in spite of the fact that just a control on the $L^2$ norm can continue the solution, as no global norm bounds are currently available.
As a matter of fact, the global-in-time existence of solutions remains essentially an open problem. Despite the fact that controlling the $L^2$ norm is sufficient to extend the solution, no global norm bounds are currently available.
The presence of the bi-harmonic operator precludes a maximum principle and energy estimates do not close to apply Gr\"onwall's inequality as in Burger's equation, since $\bm{v}$ is not divergence free. Global existence for well-prepared (small) data can be obtained if there are no growing modes, i.e., the periods are less than $2\pi$ \cite{AM19}. Otherwise, it has been shown essentially only for thin or highly anisotropic domains \cite{KukavicaKSE}, or when there is only one growing mode in each direction (a result due to one of the authors and Ambrose \cite{AM21}).  We refer the reader to \cite{LariosKSE} for a recent overview of the existing literature. Given the nature of the interplay between small and large scales in KSE, the action of a strong background flow can enhance the transfer of energy to small scale and effectively prolong the life of the solution (cf. \cite{KZquenching}). We introduce therefore the {\em advective Kuramoto-Sivashinsky equation}, which has been used as a simplified model of flame propagation in pre-mixed combustion (see for example \cite{AKSE2} and references therein):
\begin{equation} \label{eq:AKSEscalar}
  \partial_t \phi + \bm{u}\cdot \varphi + \kappa\, \big( (\Delta^2+\Delta) \varphi +|\nabla \varphi|^2 \big)=0,
\end{equation}
where $\bm{u}$ is a {\em given} vector field and the coefficient $\kappa>0$ measures the strength of the advection (in fact, by a change of time $\kappa=1/A$ where $A$ is the amplitude of the background flow).

To showcase the effect of enhanced dissipation, we present a global existence result by some of the authors for \eqref{eq:AKSEscalar}, when $\bm{u}$ is an autonomous (horizontal) shear flow with finitely many critical points of finite order $m \geq 2$, provided $\kappa$ is sufficiently small. 

Global existence then holds for arbitrary large data $\varphi_0\in L^2(\mathbb{T}^2)$. While an arbitrary number of growing modes are allowed in the direction of the shear, no growing modes can be allowed in the direction orthogonal to the shear, i.e., $L_2<2\pi$, since the advection operator has no effect on functions of $x_2$ alone. 
We stress that, while the non-linearity is the reason why global existence cannot be readily established, it is also crucial in coupling modes and, hence, propagating the effect of enhanced dissipation to all modes, even those in the kernel of the advection operator. Indeed, purely vertical modes satisfy a modified 1D KSE, which is coupled to horizontal modes only through terms that are damped by the enhanced dissipation on long-time scales. 

We recall the following result from \cite{CZDFM21}.
%Again, for notational convenience, we denote a point $\bm{x}\in \mathbb{T}^2$ as $\bm{x}=(x,y)$

\begin{theorem}[Theorem 1.1 \cite{CZDFM21}] \label{thm:KSEshear}
Let $0<L_2<2\pi$, and let $\bm{u}(x_1,x_2)=(u(x_2),0)$, where  $u:[0,L_2)\to \mathbb{R}$ is a smooth function with finitely many critical points of order at most $m\geq 2$. Let $\varphi_0 \in L^2(\mathbb{T}^2$. Then, there exists $0<\kappa_0$, which depends on $L_i$, $i=1,2$, and on $\|\varphi_0\|_{L^2}$, such that for any $0<\kappa<\kappa_0$, the AKSE \eqref{eq:AKSEscalar} admits a unique solution $\varphi \in L^\infty([0,\infty), L^2(\mathbb{T}^2))\cap L^2([0,\infty),H^2(\mathbb{T}^2))$.
\end{theorem}

Before we give a sketch of the proof, we make a few remarks. The solution is given in the class of weak solutions, but in fact any solution of KSE or AKSE  (if $\bm{u}$ is smooth enough)  with initial data in $L^2$ is a strong solution on any finite time interval. 
As already mentioned in Section \ref{s:resolvent}, the condition on the number and order of the critical points of $u$ ensures that the operator $e^{t u\partial_x}$ is mixing on the complement of its kernel, the space $\mathring{L}^2(\mathbb{T}^2)$ introduced in Section \ref{s:resolvent}. Then enhanced dissipation holds (see Theorem \ref{thm:enhancedShearTorus} and \cite{BedrossianCotiZelati17} for more details), that is, estimate \eqref{eq:sembound} holds for $L_\kappa=u\partial_{x_1} -\kappa (\Delta^2+\Delta)$ with $\lambda_\kappa=\epsilon_0\, \kappa^{2m/(2m+1)}$. Since $\lambda_\kappa/\kappa \to 0 $ as $\kappa\to 0$, this bound is sufficient to establish Theorem \ref{thm:KSEshear}. However, $\lambda'_\kappa$ in Theorem \ref{thm:ResolventEstimate}, under Assumption \ref{a:flow} on $u$,  is smaller, giving a sharper rate of decay and  a larger value for $\kappa_0$.
Lastly, we explicitly note that the decay rate in \eqref{eq:sembound} allows to estimate the dissipation time $t_{dis}$ of the system, which is discussed in Section \ref{s:CHE} (see equation \eqref{eq:dissipationTime}).

\smallskip

\noindent{\em Sketch of proof. } 
We outline briefly the main steps in the proof of Theorem \ref{thm:KSEshear}. It is based on a bootstrap argument, inspired by \cite{BH18}. 

First, by results in e.g. \cite{FM22}, we have short-time existence of a mild solution, which solves the Volterra integral formulation of the equation. The mild solution is also a weak solution and satisfies the energy identity, obtained formally by multiplying \eqref{eq:AKSEscalar} by $\varphi$ and integrating by parts. The proof of short-time existence implies that the $L^2$-norm is a continuation norm for the solution. Therefore, it is enough to obtain a global {\em a priori}
bound on $\|\varphi\|_{L^2}$.

This bound will follow from enhanced dissipation. To this effect, the solution $\varphi$ is decomposed into two parts, 
$\varphi=\varphi_\neq + \langle \varphi\rangle$, with \ $\langle \varphi\rangle(t,x_2) =\int_{0}^{L_1} \varphi(t, x_1, x_2) \, dx_1$ the projection of $\varphi$ onto the kernel of the advection operator. As in \cite{CZDFM21}, we will call $\varphi_\neq$ and $\langle \varphi\rangle$ the projected and kernel components of $\varphi$, respectively.
These two components are  coupled, but only through $\psi(x_2, t):= \partial_{x_2} \langle \psi\rangle (x_2, t)$.

Because of the local well-posedness, there exists a minimum time $t_0>0$ such that:
\begin{enumerate} [label=(H\arabic*), ref=(H\arabic*)]
\item \label{i:bootstrap1} $\|\varphi_\neq(t)\|_{L^2}\leq 8 e^{-\lambda_\nu (t-s)/4}\|\varphi_\neq(s)\|_{L^2}$,
\item  \label{i:bootstrap2} $\nu \int_s^t \|\Delta\varphi_\neq(\tau)\|_{L^2}^2\,d \tau\leq  4 \|\varphi_\neq(s)\|_{L^2}^2$,
\end{enumerate}
for all $0<s<t\leq t_0$. These are the {\em bootstrap assumptions}. 

We observe that $\psi$ satisfies a 1D KSE-like equation forced by the term $\frac{\kappa}{L_1} \int^{L_1}_0 \partial_{x_2} |\nabla \varphi_\neq|^2 \, dx$.
Hence,  the bootstrap assumptions allow to obtain a bound on $\psi$ in $L^\infty([0,t_0), L^2(\mathbb{T}))\cap L^2([0,t_0),H^1(\mathbb{T}))$, depending only on the initial data and the periods. 

This bound  and the bootstrap assumptions, in turn, give what we call the {\em bootstrap estimates}:
\begin{enumerate} [label=(B\arabic*), ref=(B\arabic*)]
\item \label{i:bootstrap3} $\|\varphi_\neq(t)\|_{L^2}\leq 4 e^{-\lambda_\nu (t-s)/4}\|\varphi_\neq(s)\|_{L^2}$,
\item  \label{i:bootstrap4} $\nu \int_s^t \|\Delta\varphi_\neq(\tau)\|_{L^2}^2\,d \tau\leq  2\|\varphi_\neq(s)\|_{L^2}^2$\,.
\end{enumerate}
To establish the bootstrap estimates, enhanced dissipation is used if $t_0$ is small, that is, when $t_0\ll t_{dis}$, the dissipation time.
Since the bounds \ref{i:bootstrap3}-\ref{i:bootstrap4} are stronger than \ref{i:bootstrap1}-\ref{i:bootstrap2} by the minimality assumption on $t_0$ it follows that  $t_0=\infty$. Therefore, we immediately have $\varphi_\neq\in L^\infty([0,\infty);L^2(\mathbb{T}^2))\cap L^2([0,\infty);H^2(\mathbb{T}^2))$, from which it also follows that $\psi\in L^\infty([0,\infty);L^2(\mathbb{T}^1))\cap L^2([0,\infty);H^2(\mathbb{T}^1))$. 
Finally, applying standard inequalities and Gr\"onwall's Lemma to the energy identity for $\varphi$ allows to conclude that additionally 
$\langle \varphi\rangle \in L^\infty([0,\infty);L^2(\mathbb{T}^1))\cap L^2([0,\infty);H^2(\mathbb{T}^1))$.
\qed

We close this section by mentioning a related global existence result when $\bm{u}$ is a {\em mixing} flow \cite{FM22}. Then there are no restrictions on the number of growing modes, that is, the size of the periodic box, since the advection operator acts on all (non-trivial) modes. The proof is simpler than in the shear flow case and it is based on a dichotomy argument. If dissipation alone is sufficiently strong, the  bound on the $L^2$ norm of the solution follows by standard energy estimates. If dissipation is not strong enough, then a control on the $L^2$ norm is achieved provided the dissipation time $t_{dis}$ is small enough, which is obtained by taking $\kappa$ sufficiently small. In fact, the global existence result does not require $\bm{u}$ to mixing explicitly, only that $t_{dis}$ is sufficiently small.

\begin{theorem}[ Main Theorem 2 \cite{FM22} ] \label{thm:KSEmixing}
Let $\varphi_0\in L^2(\mathbb{T}^2)$ and let  $\bm{u}\in L^\infty([0,\infty), L^2(\mathbb{T}^2)$. There exists $\tau_0$, depending on $L_i$, $i=1,2$, and on $\|\varphi_0\|_{L^2}$, such that, if the dissipation time $t_{dis}<\tau_0$, then there exists a unique mild and weak solution  $\varphi \in  L^\infty([0,T);L^2(\mathbb{T}^2))\cap L^2([0,T);H^2(\mathbb{T}^2))$ of \eqref{eq:AKSEscalar} for any $T>0$.
\end{theorem}

\printbibliography
\end{document}